\documentclass[12pt]{iopart}

\usepackage{iopams}  
\expandafter\let\csname equation*\endcsname\relax
\expandafter\let\csname endequation*\endcsname\relax
\usepackage{amsmath}
\usepackage{amssymb}
\usepackage{graphicx}
\usepackage{epstopdf}
\usepackage{appendix}
\usepackage{stmaryrd}
\usepackage{url}
\usepackage{hyperref}

\usepackage[utf8]{inputenc}
\usepackage{alphabeta}

\begin{document}

\title[Discontinuous collocation methods \& gravitational self-force applications]{Discontinuous collocation methods and gravitational self-force applications}

\author{Charalampos Markakis$^{1,2,3}$, Michael F. O'Boyle$^{4}$, Pablo D. Brubeck$^{5}$, Leor Barack$^6$}

\address{$^{1}$DAMTP, Wilberforce Road, University of Cambridge, Cambridge, CB3 0WA, UK\\
        $^{2}$School of Mathematical Sciences, Queen Mary University of London, Mile End Rd, London, E1 4NS, UK\\
        $^{3}$NCSA, University of Illinois at Urbana-Champaign,
        1205 W Clark St, Urbana, Ilinois 61801, USA\\
        $^{4}$Department of Physics, University of Illinois at Urbana-Champaign, 1110 W Green St, Urbana, IL 61801, USA\\
        $^{5}$Mathematical Institute, University of Oxford, Oxford, OX26GG, UK \\
        $^{6}$Mathematical Sciences, University of Southampton,
        Southampton, SO17 1BJ, UK\\
        }
\ead{$^1$c.markakis@damtp.cam.ac.uk,$^4$moboyle2@illinois.edu,
$^5$brubeckmarti@maths.ox.ac.uk,
$^6$l.barack@soton.ac.uk}
\vspace{10pt}
\begin{indented}
\item[]December 2020
\end{indented}

\begin{abstract}
Numerical simulations of extereme mass ratio inspirals, the mostimportant sources for the LISA detector, face several computational challenges. We present a new approach to evolving partial differential equations occurring in black hole perturbation theory and calculations of the self-force acting on point particles orbiting supermassive black holes. Such equations are distributionally sourced, and standard numerical methods, such as finite-difference or spectral methods, face difficulties associated with approximating discontinuous functions. However, in the self-force problem we typically have access tofull a-priori information about the local structure of the discontinuity at the particle. Using this information, we show that high-order accuracy can be recovered by adding to the Lagrange interpolation formula a linear combination of certain jump amplitudes. We construct discontinuous spatial and temporal discretizations by operating on the corrected Lagrange formula. In a method-of-lines framework, this provides a simple and efficient method of solving time-dependent partial differential equations, without loss of accuracy near moving singularities or discontinuities. This method is well-suited for the problem of time-domain reconstruction of the metric perturbation via the Teukolsky or Regge-Wheeler-Zerilli formalisms. Parallel implementations on modern CPU and GPU architectures are discussed.
\end{abstract}

%
%
%
%
%

\section{Introduction}
In black hole perturbation theory, the metric tensor perturbation $\delta g_{\mu \nu}$
can be reconstructed from a  scalar master function which satisfies the Teukolsky equation \cite{Teukolsky1973}. For a Schwarzschild spacetime background, upon (tensor) spherical harmonic decomposition, in geometric ($G=c=1$) units, this equation can be written as a  sourced  wave equation:
\begin{equation} \label{TeukSchw}
\negthickspace[ - \partial _t^2 + \partial _x^2 + {A^t}(x){\partial _t} + {A^x}(x){\partial _x} + {V_l}(x)]{\Psi _{lm}}(t,x) 
= {G_{lm}}(t,x)\delta (x - \xi (t)) + {F_{lm}}(t,x)\delta '(x - \xi (t)),
\end{equation}
where $x=r+2M \ln(\frac{1}{2}r/M-1)$ is the tortoise coordinate, $t,r$ are  Schwarzschild coordinates and $x=\xi(t)$ is the time-dependent radial location of the particle. The Regge-Wheeler \cite{Regge1957} and  Zerilli \cite{Zerilli1970a,Zerilli1970} equations have a similar form, with $A^\mu$ set to zero. A similar equation governs perturbations in Kerr spacetime, albeit with an additional source term $\delta''(x - \xi (t))$ and coupling between $l$-modes. The latter issue can be dealt with the use of Clebsch-Gordan coefficients and inverting/storing the mode-coupling matrix \cite{Dolan2013}.

In calculations of radiation-reaction forces  in electrodynamics, or gravitational
self-force corrections to geodesic motion of small bodies in general relativity
   \cite{Barack2009,Poisson2011},  the  source of the   dimensionally reduced
field equations  is singular, making the solution or its derivatives discontinuous
across the particle location. Methods of solving distributionally-sourced
equations
for  self-force applications  have hitherto included  \textit{(i)} multi-domain pseudospectral methods  
with time-dependent coordinate mapping   and junction conditions at the discontinuity
\cite{Canizares2009,Canizares2011a,Canizares2010c,Jaramillo2011}, \textit{(ii)} discontinuous
Galerkin methods with time-dependent coordinate mapping   and junction conditions at the discontinuity \cite{Field2009,Field:2010xn,Heffernan:2017cad},
\textit{(iii)} finite-difference methods based on null coordinates
which incorporate the source term without approximating the Dirac $\delta$ or $\delta'$ functions
 \cite{Barack:2000zq,Lousto:1997wf,Lousto:1999za,Barack:2002ku,Barack:2005nr,Haas:2007kz,Haas:2011np,Martel:2001yf}
\textit{(iv)}
finite-difference representations of the Dirac  $\delta$ function 
 \cite{Sundararajan2007,Sundararajan2008}
and \textit{(v)}
Gaussian
function representations of the Dirac  $\delta$ function  \cite{Sundararajan2007,Sundararajan2008,Harms2013}.
Methods (i)-(ii) are based on a Lagrangian picture with coordinates that track the particle, while methods (iii)-(v) are based on an Eulerian picture, with fixed coordinates and the particle moving freely in the computational grid. It would be highly desirable to construct a scheme with the ease and flexibility of Eulerian methods and the high precision of Lagrangian methods. The discontinuous collocation methods constructed below achieve  this by extending and systematizing earlier approaches. In particular,  it will be shown that inclusion of known derivative jumps of arbitrarily high order allows finite-difference or pseudospectral implementations in Eulerian coordinates without domain decomposition.

Characteristic finite-difference methods have been highly successful, but their extension, for instance, from second to fourth order involves considerable effort and  extensive code rewriting. 
Finite-difference and spectral methods,  commonly used for numerically solving ordinary or partial differential equations, are highly efficient and accurate when the solution is
smooth. When discontinuities are present, however,  oscillations  due to  Gibbs phenomena arise that contaminate the solution  and cause loss of  point-wise convergence. Several  workarounds  exist and are commonly used to suppress or avoid these phenomena. However, they typically come at a cost. For piecewise analytic solutions, domain decomposition is highly accurate as it avoids the point of discontinuity and  recovers the lost convergence, but numerical implementation  can be complicated:   for certain numerical-relativity problems that involve moving discontinuities, the requirements of remapping of coordinates or regridding  at each time step add to the computational cost and complexity.  For finite-difference methods,  lowering the order of  approximation globally or near the discontinuity  reduces oscillations but sacrifices accuracy. For spectral methods, applying filtering techniques can improve the  convergence properties of the solution away
from discontinuities, but the accuracy near  discontinuities remains poor \cite{Boyd2001}. Post-processing of the
resulting oscillatory data using direct \cite{Gottlieb1992,Gottlieb1994,Gottlieb1996,Gottlieb1995,Gottlieb1995a,Gottlieb1997,GottliebDavidSigal2003,2019JCoPh.390..527P} and inverse \cite{Jung2004} Gegenbauer reconstruction
can  recover spectrally convergent non-oscillatory solutions. This approach has been successfully used in  several contexts, such as magnetic resonance imaging
\cite{Archibald2003,Archibald2002,Archibald2004}
and evolving nonlinear hyperbolic systems \cite{Shu1995}.  Nevertheless, reconstruction techniques are   computationally
expensive (the cost of reconstruction is $O(N^2)$  where $N$ is the
number of gridpoints) and  subject to  complications. For instance, the Gegenbauer transformation matrices become ill-conditioned for high-degree polynomials \cite{Jung2004}.
In addition, the 
 reconstruction must  be constrained to use a sufficiently small ratio of order  to truncation, which can  impair its effectiveness  in certain applications \cite{Boyd2005}.

Methods of approximating piecewise smooth functions  have been introduced by  
Krylov \cite{Krylov1907}, Lanczos \cite{Lanczos1966} and  
Eckhoff \cite{eckhoff1994,eckhoff1995,eckhoff1996,eckhoff1997,eckhoff1998,eckhoff2002}
  in the context of Fourier methods (see also \cite{kantorovichkrylov1958,Barkhudaryan2007,Poghosyan2010,Adcock2010,Batenkov2012,Poghosyan2012})
and by Lipman and Levin \cite{Lipman2009} in the context of finite-difference
methods. These methods have used a polynomial correction term, written as an  expansion in a monomial basis, with  coefficients determined from suitable jump conditions. Lipman and Levin moreover introduced a moving least squares method for determining the location and magnitude of the discontunuity which scales as $O(N)$, rather than $O(N^2)$ as with Gegenbauer reconstruction methods.

For   self-force problems, however, the location and  magnitude of the discontinuity can be computed \textit{a priori} at no cost.
Earlier approaches have used jumps of the function and its low order derivatives. In this paper, we extend and systematize the use of jumps to arbitrary derivative order.
In particular, as shown in Appendix A, jumps of the master function and all its derivatives   
  can be determined exactly from $G_{lm}$ and $F_{lm}$.
We construct a scheme particularly suited  for numerically solving such problems, which exploits the fact that the location and magnitude of the jumps are known. 

In what follows, we focus on  nodal  methods based on Lagrange interpolation. We show that adding a linear combination of the jumps to the Lagrange interpolation formula can account for the presence of the discontinuity and restore the original accuracy. Our work  extends  earlier piecewise interpolation formulae  by introducing a Lagrange basis, rather than a monomial basis. This is vital, as it allows one to go beyond interpolation, and to construct  differentiation matrices or quadrature weights  for discontinuous collocation methods. 
Finite-difference and pseudo-spectral scheme developed for smooth problems, outlined in \S2,  can thus be  straightforwardly extended  to discontinuous  problems, illustrated for Schwarzschild perturbations in \S3. 

In a method-of-lines context, the discontinuous collocation methods described above are used for spatial discretization. However, when used with standard time integration methods developed for smooth problems, such as Runge-Kutta methods, this approch works only for static particle sources. For moving discontinuities, a discontinuous time integration method is necessary. In \S4, we thus  develop a discontinuous time integration formula that takes into account known jumps accross the particle world-line. The formula presented here in is 2$^\text{nd}$ order in time, and can be considered a discontinuous generalization of the trapezium rule. 4$^\text{th}$ and higher order formulas can also been constructed, as we will demonstrate in a companion paper. 
Merits, possible extensions and foreseeable applications of this method are summarized in \S5.

\section{Smooth functions}
As a primer for the discontinuous collocation methods that follow in \S3, we first review the Lagrange interpolation formula  for sufficiently differentiable functions and its role in numerical integration and differentiation.

\subsection{Interpolation}
Let $f:[a,b]\rightarrow \mathbb{R}$  be a $C^{N+1}$ function, with its values $f_i=f(x_i)$ known at the ordered, distinct nodes  $x_i$
($a\le x_0<x_1<\dots<x_N \le b$). 
The collocation  polynomial of degree $N$,
\begin{equation} \label{PolyN1}
p(x) = \sum\limits_{j = 0}^{N} {{c_j}{x^j}} ,
\end{equation}
which interpolates the given $N+1$
nodes, is unique and determined by solving the linear system of algebraic conditions 
\begin{equation} \label{PolyNConds1}
p(x_{i}) = f_i 
\end{equation}
for the coefficients $c_j$. The solution  may be written in the Lagrange form
\begin{equation} \label{PolyNLagrange1}
p(x) = \sum\limits_{j = 0}^N {{f_j}{\pi _j}(x),} 
\end{equation}
where
\begin{equation} \label{PolyNLagrangeBasis1}
{\pi _j}(x) = \prod\limits_{\genfrac{}{}{0pt}{2}{k = 0}{
k \ne j}}^N {\frac{{x - {x_k}}}{{{x_j} - {x_k}}}} 
\end{equation}
are the Lagrange basis polynomials. By construction, the basis polynomials satisfy
the conditions $\pi_j(x_i)=\delta_{ij}$, with $\delta_{ij}$ denoting the Kronecker symbol, so that the  polynomial \eqref{PolyNLagrange1} satisfies the conditions \eqref{PolyNConds1}.

The nodes have so far been left unspecified. For polynomials of low degree, equidistant nodes
\begin{equation} \label{EquidistantNodes}
x_i=a+i\frac{{b- {a}}}{N}, \quad i=0,1,\dots,N
\end{equation}
are  commonly used. For polynomials of high degree, however, this choice leads to  oscillations due to  Runge's phenomenon and the interpolation process may diverge. 

The oscillations can be minimized by using nodes that are distributed more densely towards the edges of the interval. For instance, interpolation on Chebyshev-Gauss-Lobatto nodes
\begin{equation}  \label{ChebyshevNodes}
x_i=\frac{{{a} + {b}}}{2} + \frac{{{a} - {b}}}{2}\cos \theta_i, 
\quad
\theta_i= {\frac{{i\pi }}{N}}, \quad i=0,1,\dots,N
\end{equation}
converges uniformly 
for every absolutely continuous function \cite{Fornberg1998}. These  nodes are the extrema of a Chebyshev polynomial of degree $N$. If all the function derivatives are bounded on the  interval $x \in [a,b]$,
 Chebyshev interpolation always has the property of infinite order (i.e. exponential)  convergence. 

Spectral methods are often introduced and implemented  as 
 expansions in terms of orthogonal basis functions.
However, in many practical applications one introduces a finite grid, which requires evaluation of basis functions on each grid point. Nevertheless, it is often unnecessary to explicitly compute the basis functions. Indeed, the  construction of a spectral expansion, with  coefficients determined via collocation at the nodes  $x_i$, can be replaced by  a global interpolating polynomial if the nodes are the abscissas of a  Gaussian
quadrature (typically the simple roots or extrema of a Jacobi  polynomial, of which Chebyshev is a special case)
associated with the basis functions \cite{Fornberg1998,Boyd2001}. This pseudospectral approach   allows numerical operations to be performed  in physical, rather than spectral, space without compromise in convergence order. This property can be exploited in     numerical codes, offering the flexibility of switching between finite-difference and pseudospectral methods   by simply shifting the location of the 
nodes [e.g. changing from Eq.~\eqref{EquidistantNodes} to \eqref{ChebyshevNodes}]. We shall make use of this flexibility in \S3.

%




\subsection{Differentiation}
The $n$-th derivative $f^{(n)}=d^n f/dx^n$ of $f$
may be approximated by differentiating the Lagrange interpolating polynomial 
\eqref{PolyNLagrange1} $n$ times. This operation may be written as matrix multiplication:
\begin{equation} \label{fnderivativeapprox1}
{f^{(n)}}({x_i}) \simeq {p^{(n)}}({x_i }) = \sum\limits_{j = 0}^N D_{ij}^{(n)}f_j,
\end{equation}
where 
\begin{equation}  \label{derivativematricespi}
D_{ij}^{(n)} = \pi _j^{(n)}({x_i})
\end{equation}
is the $n^{\rm{th}}$ order \textit{differentiation matrix}, constructed from the 
$n^{\rm{th}}$ derivative of the basis polynomials
\eqref{PolyNLagrangeBasis1} 
with respect to $x$ at $x_i$. Note that, for $n=0$, the above expression gives the identity matrix, $D_{ij}^{(0)}=\delta_{ij}$,  as it ought to.
For fixed grids, the above matrices can be pre-computed and stored in memory.  

Fast methods for  computing the  matrices 
 have been developed
\cite{Fornberg1998,Fornberg1988,Fornberg2006,Sadiq2011,Welfert1997} and implemented in  computational libraries and software. For instance, the Wolfram Language command
\begin{verbatim}
Dn = NDSolve`FiniteDifferenceDerivative[Derivative[n],X,
"DifferenceOrder" -> m] @ "DifferentiationMatrix"
\end{verbatim}
uses Fornberg's algorithm \cite{Fornberg1998,Fornberg1988,Fornberg2006} to compute the  
$n$-th order differentiation matrix \verb Dn ,
given the list of grid-points 
$X =\{x_0,x_1,x_2,...,x_N\}$,   at
the desired order of approximation 
 $m$. The latter is an integer satisfying $n \le  m \le N$ that yields a finite-difference error $\mathcal{O}(\Delta x ^ { m } 
)$, where $\Delta  x $ is the maximum local grid spacing.

A choice $m<N$ yields a   composite derivative matrix of size $N+1$ constructed from centred 
$(m+1)$-point stencils. For example, the choice $n=1$ and $m=2$ computes $1^{\rm st}$ derivatives of  $2^{\rm nd}$-order polynomials
[given by Eq.~\eqref{PolyNLagrange1} for $N=2$] to 
 construct the centred finite-difference approximation
\begin{eqnarray}  \label{smoothFDf}
{{f'}_i} &=& {f_{i + 1}}\frac{{{x_i} - {x_{i - 1}}}}{{({x_{i + 1}} - {x_i})({x_{i + 1}} - {x_{i - 1}})}} + {f_i}\frac{{{x_{i + 1}} + {x_{i - 1}} - 2{x_i}}}{{({x_{i + 1}} - {x_i})({x_i} - {x_{i - 1}})}} 
\\
&-& {f_{i - 1}}\frac{{{x_{i + 1}} - {x_i}}}{{({x_{i + 1}} - {x_{i - 1}})({x_i} - {x_{i - 1}})}}+\mathcal{O}(\Delta x^2).  \nonumber
\end{eqnarray}
The entire domain $[a,b]$ is covered in a composite fashion (that is, the above formula is applied repeatedly at each node $x_i$ using only nearest-neighbour points). This gives the centred-derivative matrix
\begin{eqnarray}
D_{ij}
&=&\frac{{{x_i} - {x_{i - 1}}}}{{({x_{i + 1}} - {x_i})({x_{i + 1}} - {x_{i - 1}})}}{\delta _{j,i + 1}} + \frac{{{x_{i + 1}} + {x_{i - 1}} - 2{x_i}}}{{({x_{i + 1}} - {x_i})
({x_i} - {x_{i - 1}})}}{\delta _{ij}} 
\\
&-& 
\frac{{{x_{i + 1}} - {x_i}}}{{({x_{i + 1}} - {x_{i - 1}})({x_i} - {x_{i - 1}})}}{\delta _{j,i - 1}} \nonumber  
\end{eqnarray}
 for all $i=1,\dots,N-1$, while one-sided finite differences are used  at the end-points $i=0,N$.

A single-domain pseudospectral derivative matrix can be obtained when the grid-points $x_i$ are spaced corresponding to  the Chebyshev-Gauss-Lobatto nodes \eqref{ChebyshevNodes}. Then, the elements  $D_{ij}$ of  the $(N+1)\times (N+1)$ Chebyshev derivative
matrix $D$ are given by
\begin{eqnarray}
D_{00} &= & -\frac{2N^2+1}{3(b-a)}, \quad i=j=0 \nonumber \\
D_{i i} &=&  \frac{\cos \theta_i}{(b-a)\sin^2 \theta_i}, \quad 1\le i \le N-1   \nonumber \\
D_{NN} &= & \frac{2N^2+1}{3(b-a)}, \quad i=j=N \nonumber \\
D_{i j} &= & \frac{2c_i (-1)^{i+j}}{(b-a)c_j (\cos \theta_j-\cos \theta_i)}, \quad 0\le i,j \le N \,\,{\rm{and}}
\,\, i\ne j \nonumber 
\end{eqnarray}
where $\theta_i = i \pi/N$ and $c_i=1+\delta_{i,0}+\delta_{i,N}$ (that is, $c_i=1$ for $1 \le i \le
N-1$ and $c_0=c_N=2$). Round-off error of numerical differentiation with the above matrix can be reduced via the  negative sum trick of Ref.~\cite{Baltensperger2003},
which amounts to replacing the  diagonal elements  $D_{ii}$ by $- \sum_{j\ne i} D_{ij}$. This ensures that,   when  the derivative
matrix multiplies a constant
vector, the result is exactly zero, that is, $\sum_j D_{ij}=0$.

The Chebyshev derivative matrix can be constructed automatically via the Wolfram Language command given above by using the maximal value 
$m=N$  when computing the differentiation matrix, or, equivalently, using
the command:
 \begin{verbatim}
Dn = NDSolve`FiniteDifferenceDerivative[Derivative[n],X,
"DifferenceOrder" -> "Pseudospectral"] @ "DifferentiationMatrix"
\end{verbatim}
which returns the order-$n$ derivative matrix $D^n$. 

The pseudospectral derivative matrices obtained from this (single-domain)  application of Eq.~\eqref{derivativematricespi}  are the \textit{infinite-order limit} of the finite-difference derivative matrices obtained by  (composite) application of   Eq.~\eqref{derivativematricespi} 
\cite{Fornberg1998}.
In view of the fact that pseudo-spectral methods are a limiting case of finite-difference methods, we find it advantageous to construct our code using differentiation matrices. This allows one the flexibility of specifying the grid and desired order of the method at the beginning of the code, while the rest of the code remains the same, and works equally well for finite-difference or pseudo-spectral methods.
  

\subsection{Integration}
Quadrature, or numerical integration, amounts to approximating the definite integral of a given function or, equivalently,  integrating an approximation to the function. Integrating the approximation \eqref{PolyNLagrange1} gives a quadrature rule in the form of a weighted sum of function evaluations:
\begin{equation} \label{Integralp1}
\int_{{a}}^{{b}} {f(x)dx}\simeq  \int_{{a}}^{{b}} {p(x)dx}  = \sum\limits_{j = 0}^N {{w_{j}}{f_j}},
\end{equation}
where
\begin{equation} \label{IntegralMatrixp1}
{w_{j}}=\int_{{a}}^{{b}} {\pi_j(x)dx}  
\end{equation}
are the relevant weights of integration.
Various well-known quadrature rules stem from the above formulae. For instance, using interpolation of 2 or 3  equidistant nodes  yields  
the composite trapezoidal or Simpson rule respectively.
For the Chebyshev-Gauss-Lobatto nodes \eqref{ChebyshevNodes}, the above procedure  yields the homonymous quadrature formula, which is exponentialy convergent for smooth integrands. 

\section{Non-smooth functions} 
In general, a jump discontinuity in a  function or its derivatives
cannot  be  accurately approximated by a single  polynomial.
  As mentioned earlier, a common workaround is to divide the computational domain so that  any discontinuities  are located at the boundary between domains. 
In what follows, we explore a simpler approach that does not require domain decomposition  but uses the location and magnitude of the discontinuity if known.
We will illustrate that incorporating discontinuities into the approximating function instead of avoiding them offers flexibility and ease of  programming without loss of accuracy.
We will use a `nodal' approach, and the method of undetermined coefficients, to accommodate the necessary collocation and jump conditions. An equivalent `modal' approach is described in Appendix A.

\subsection{Interpolation} \label{InterpolationSection}
Let  $f:[a,b]\rightarrow \mathbb{R}$  now be a $C^{k}$ function,  with the values $f_i=f(x_i)$ known at the ordered, distinct nodes  $x_i$.  The Lagrange interpolation method outlined in \S2 is limited by the degree of differentiability, that is,    $N+1 \le k$.
Any  degree of differentiability $k$ lower than $N+1$ necessitates  modification   of the above methods. 
If  $f$ is a piecewise-$C^{N+1}$ function and  its jump discontinuities are known, then a simple generalization of the Lagrange method to
 any   
$k \ge -1$ may be obtained as follows. 

Let  the discontinuity be located at some $\xi \in (a,b)$ and  the (non-infinite) jumps in $f$ and its derivatives  
 be given by:
\begin{equation} \label{fmJconditions2}
{f^{(m)}}(\xi^+  ) - {f^{(m)}}(\xi^- ) = {J_m},\quad m = 0,1, \ldots ,\infty,
\end{equation}
where
$\xi^-$ and $\xi^+$ denote the left and right limits to $\xi$.
One may then approximate $f$ 
by a piecewise polynomial function 
that (i) interpolates the given $N+1$ nodes, and (ii)  has the same derivative jumps at the discontinuity. Eq. \eqref{PolyN1}  is thus replaced by the ansatz
\begin{equation} \label{PolyN2DiscInterpol}
p(x) 
=\theta (x - \xi ){p_ + }(x) + \theta (\xi  - x){p_ - }(x)
\end{equation}
where
\begin{equation} \label{PolyN1DiscInterpol}
\theta (x)=\left\{ {\begin{array}{*{20}{c}}
1,&{x > 0}\\
{1/2,}&{x = 0}\\
0,&{x < 0}
\end{array}} \right.
\end{equation}
denotes the Heaviside step function,
and 
\begin{equation} \label{PolyN2}
{p_ - }(x) = \sum\limits_{j = 0}^{N} {c_j^-(\xi) {x^j}},\quad{p_ + }(x) = \sum\limits_{j = 0}^{N} {c_j^+(\xi) {x^j}} 
\end{equation}
denote the left and right interpolating polynomials.
Note that the polynomials $p_-,p_+$  and the piecewise polynomial function $p$ depend implicitly  on the location $\xi$ of the discontinuity.

Half of the ($2N+2$)   polynomial coefficients $c_j^\pm$
are to be determined by the collocation conditions
\eqref{PolyNConds1} which, given the ansatz \eqref{PolyN2DiscInterpol}, translate to
\begin{equation} \label{PolyNConds2}
\left\{ \begin{array}{*{20}{c}}
p_{+} (x_i) = f_i,&x_i > \xi\\
p_{-} (x_i) = f_i,&x_i < \xi
\end{array} \right.
.
\end{equation} 
The remaining  coefficients could be determined by the first $N+1$ of the jump conditions  \eqref{fmJconditions2}. However, as information on  high order jumps may not be   beneficial (as discussed below and illustrated in Fig.~3) or easily accessible,  we shall allow the number of  jumps enforced  to be user-specifiable, by dropping jumps in derivatives higher than order $M\le N$, where $M$ is some chosen integer between
$-1$ and $N$. 
 The remaining $N+1$ coefficients are  then determined by the jump conditions
\begin{equation} \label{fmJ0conditions2}
p_ + ^{(m)}(\xi ) - p_ - ^{(m)}(\xi ) = 
\left\{ \begin{array}{l}
{J_m},\quad m = 0,1,...,M\\
0,\quad m = M + 1,...,N
\end{array} \right.
.
\end{equation}
Setting the higher than $M^{\rm th}$-order  derivative jumps  to zero  means that the polynomial coefficients $c_{M+1},\dots, c_N$ remain constant across $\xi$ as in  Lagrange interpolation. This choice is favoured by numerical experiments to be discussed below and allows recovery of Lagrange interpolation as a special case, $M=-1$. 

It should be noted that the ansatz \eqref{PolyN2} should not be viewed merely as two independent interpolating polynomials matched at the discontinuity. If that were the case, the polynomial degrees would differ  left and right of the discontinuity, depending on the number of grid points on each side, and our method would amount to domain decomposition. Instead, our ansatz should be regarded as a single ``polynomial'' of degree $N$ that covers the whole domain $[a,b]$ and has piecewise constant coefficients with  a jump  at the discontinuity 
$\xi$. The left and right coefficients encode information from \textit{all} nodes and are correlated via the jump conditions \eqref{fmJ0conditions2}.

 As mentioned above, the ``no jump'' choice $M=-1$  corresponds to 
Lagrange interpolation, which   suffices for smooth functions.  
A choice of $M\ge k+1 \ge 0$  is required to approximate  $C^{k}$ functions with  accuracy better than Lagrange,  as this reproduces the first nonvanishing jump, $J_{k+1}$, in the $(k+1)^{\rm th}$ derivative. [For instance, a $C^0$ function can be approximated by   specifying at least the first derivative jump $f'(\xi^+)-f'(\xi^-)=J_1$.]
 Incorporating higher derivative jumps, that is, using  $M > k+1 $, increases the approximation accuracy of the function  as  the  jumps $J_{k+2}, J_{k+3},\dots ,J_{M}$ in its $(k+2)^{\rm th},(k+3)^{\rm th}, \dots , M^{\rm th}$ derivative are also reproduced.\footnote{Heuristically, such behaviour is similar to Hermite \cite{hermite1878} and Birkhoff \cite{Birkhoff1906} interpolation, which can also be obtained with the method of undetermined coefficients, except the present method only uses  derivative \textit{jumps} at $\xi$ rather than function derivatives at $x_i$. A numerical convergence analysis is provided in Fig.~3, but a rigorous error analysis is beyond the scope of this paper and is left for future work.}
Nevertheless, the accuracy degrades for  values of $M$ too close to $N$ due to the appearance of  spurious  oscillations analogous to Runge's phenomenon. Our numerical experiments  suggest that intermediate choices of $M$  relative to $N$  (typically $M/N \approx 1/2 $ depending on the function) are optimal and consistently yield  highly accurate results.

As detailed above, the linear system of algebraic conditions
\eqref{PolyNConds2}-\eqref{fmJ0conditions2} can be uniquely solved by substituting the  polynomials \eqref{PolyN2} and determining their unknown coefficients.
This procedure is  simplified if  the polynomials are first converted to a Lagrange basis, 
\begin{equation} \label{ppmCpi}
p_\pm(x)=\sum_{j=0}^N C^\pm_j(\xi) \pi_j(x).
\end{equation}
Substitution in the algebraic conditions \eqref{PolyNConds2}-\eqref{fmJ0conditions2} yields
\begin{equation} \label{Cpmxifthetag}
C^+_j(\xi) =f_j+\theta(\xi - x_j)g(x_j-\xi) ,\quad C^-_j(\xi) =f_j-\theta(x_j-\xi)g(x_j-\xi),
\end{equation}
where
\begin{equation} \label{gjofx2}
g(x_j-\xi) =  \sum\limits_{m = 0}^{M } {\frac{{{J_m}}}{{m!}}{{({x_j} - \xi )}^m}} 
\end{equation}
are weights computed from the jump conditions at the discontinuity $\xi$ given the nodes $x_j$.
Eq.~(\ref{gjofx2}) can be treated as a polynomial in $\Delta \xi_j = x_j-\xi$ and computed with significantly higher speed and numerical precision in Horner form:
\begin{equation} \label{gjofx3}
g(\Delta \xi) =  J_0+\Delta \xi \left( J_1+\Delta \xi \left(\frac{J_2}{2!}+\Delta \xi \left(...+\Delta \xi\left(\frac{J_{M-1}}{(M-1)!}+\Delta \xi \frac{J_{M}}{M!}\right)...\right)\right)\right ).
\end{equation}
Substituting 
Eqs. \eqref{ppmCpi}-\eqref{gjofx2} into  the ansatz \eqref{PolyN2DiscInterpol}
yields the concise solution   
\begin{equation} \label{PolyNLagrange2}
p(x) = \sum\limits_{j = 0}^N {{[f_j+\Delta(x_j-\xi;x-\xi) ]}{\pi _j}(x)},
\end{equation}
which generalizes the Lagrange formula \eqref{PolyNLagrange1} to a \textit{piecewise} \textit{polynomial} interpolating function.
This interpolating function depends implicitly on the location $\xi$ of the discontinuity through the 2-point functions
\begin{equation} \label{sjofx2}
\Delta(x_j-\xi;x-\xi) = [\theta (x - \xi )\theta (\xi  - {x_j}) - \theta (\xi  - x)\theta ({x_j} - \xi )]g(x_j-\xi). 
\end{equation}
When  $J_m = 0$, the $\Delta$ functions vanish and  Lagrange
interpolation is recovered. When evaluated at a node $x=x_i$, the above expression simplifies to
\begin{equation} \label{sjofx2b}
\Delta(x_j-\xi;x_i-\xi)  = [\theta (x_{i} - \xi ) - \theta ({x_j} - \xi )]g(x_j-\xi).
\end{equation}
Since the prefactor in $\Delta(x_j-\xi;x_i-\xi) $ is antisymmetric in its two indices and $\pi_j(x_i)=\delta_{ij}$ is symmetric, the correction 
$\sum_j \Delta(x_j-\xi;x-\xi) \pi_j(x) $ in the interpolating formula  
\eqref{PolyNLagrange2} 
vanishes at each node $x=x_i$.
This ensures that  the collocation conditions \eqref{PolyNConds1} or, equivalently, \eqref{PolyNConds2}, are satisfied by construction.


\subsection{Static  charge in Schwarzschild spacetime} 

As an illustrative example, let us consider the  inhomogeneous Legendre differential equation with a singular (Dirac delta function) source:
\begin{equation} \label{LegendreEq}
(1-x^2)\frac{d^2\Phi_l(x)}{dx^2}-2x \frac{d\Phi_l(x)}{dx}+l(l+1)\Phi_l(x)=G_l(x)\delta(x-\xi).
\end{equation}
 Such equations arise in the context of self-forces acting on static scalar and electric test charges in the spacetime of a Schwarzschild black hole 
[cf. Ref.~\cite{Burko2000} for the derivation and the source functions $G_l(x)$]. Equation~\eqref{LegendreEq} admits 
piecewise smooth
solutions
\begin{equation} \label{PhilLegendrePQ}
\Phi_l(x)={P_l}(\xi ){Q_l}(x)\theta (x - \xi ) + {P_l}(x){Q_l}(\xi )\theta (\xi-x),
\end{equation}
where $P_l(x)$ and $Q_l(x)$ denote the Legendre functions of the first and second kind respectively.
In general, although   such solutions have finite differentiability at the test charge location $x=\xi$,  their derivative jumps can  be obtained analytically  from the field equation \eqref{LegendreEq} they satisfy without knowledge of the actual solution \cite{Burko2000,Field2009}. In this case, the derivative  jumps are given by
\begin{equation} \label{Phijumps}
J_k=\Phi^{(k)}_l(\xi^+)-\Phi^{(k)}_l(\xi^-)={P_l}(\xi )Q^{(k)}_l(\xi) - P^{(k)}_l(\xi){Q_l}(\xi ).
\end{equation}
This knowledge of the jumps can be used to  inform the method outlined above and   obtain numerical solutions to the radial part of the field equations. 
The accuracy of the numerical solutions  depends on how well the underlying interpolating function  (used for numerical operations such as integration or differentiation) approximates the actual solution. 

A comparison of the analytic solution \eqref{PhilLegendrePQ} to  
 interpolating functions based on  the Lagrange formula \eqref{PolyNLagrange1} and its generalization, Eq.  \eqref{PolyNLagrange2},
is depicted in Fig.~1 for equidistant and Fig.~2 for pseudospectral nodes.
Although both formulae exhibit Runge oscillations near the boundaries for equidistant nodes, only the corrected formula faithfully  reproduces the discontinuity given sufficient resolution. Pseudospectral nodes eliminate the Runge oscillations near the boundaries, but  the corrected  formula manages to also eliminate the Gibbs oscillations of the Lagrange formula near the discontinuity. Evidently, \textit{the accuracy  lost when applying the Lagrange formula on nonsmooth functions is  recovered  by simply  adding the correction terms} \eqref{sjofx2}.

 A comparison of convergence rates using pseudospectral nodes and different numbers of derivative jumps is depicted in Fig.~3. As mentioned earlier, the convergence estimates favour the use of an intermediate number of derivative jumps relative to a given number of nodes $(M\approx N/2)$. With pseudospectral nodes,  the  formula  \eqref{PolyNLagrange2}
is at least $M^{\rm {th}}$-order convergent  if  derivatives  $f^{(n>M)}$ have finite jumps,  as in the example given above. For functions whose derivatives higher than order $M$ are continuous, the  formula   \eqref{PolyNLagrange2} is exponentially  convergent.

\begin{figure}[ht] 
\includegraphics[width=\columnwidth]{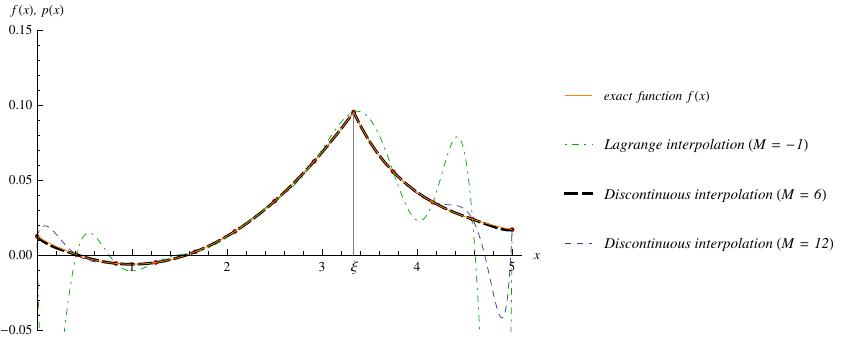}
\includegraphics[width=\columnwidth]{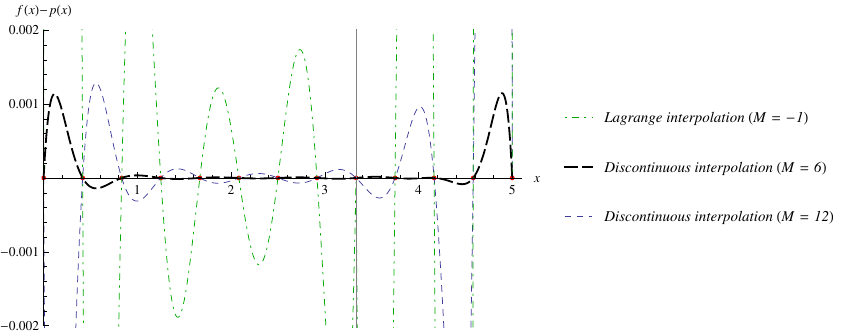}
\caption{\textbf{Top}: Single domain interpolation of the function ${f}(x) = {P_2}(\xi ){Q_2}(x)\theta (x - \xi ) + {P_2}(x){Q_2}(\xi )\theta (\xi-x)$  at $N+1=13$ equidistant nodes, using the Lagrange formula~\eqref{PolyNLagrange1} or its discontinuous generalization \eqref{PolyNLagrange2} with jump conditions on the first $M=6$ or $12$ derivatives. Although all formulas exhibit Runge oscillations near the  boundaries, only  Eq.~\eqref{PolyNLagrange2}  faithfully  reproduces the function near the discontinuity. \textbf{Bottom}: Interpolation error.  Runge oscillations are minimized when an intermediate number of jumps ($M\approx N/2$) is used.} 
\label{diff} 
\end{figure}  

\subsection{Discontinuous differentiation} A finite-difference or pseudospectral approximation to the $n$-th derivative of $f$ can be obtained by differentiating the 
piecewise polynomial \eqref{PolyNLagrange2}.
When evaluated at a node, $x=x_i$, this yields
\begin{equation} \label{fnderivativeapprox2}
{f^{(n)}}({x_i}) \simeq {p^{(n)}}({x_i }) = \sum\limits_{j = 0}^N D_{ij}^{(n)}[f_j+
\Delta(x_j-\xi;x_i-\xi)],
\end{equation}
with the differentiation matrices $D_{ij}^{(n)}$ given again by Eq.~\eqref{derivativematricespi}. As shown in Sec.~\ref{sec:sec4} below, high-order (or exponential) convergence can be attained thanks to the 2-point $\Delta$ functions  \eqref{sjofx2b}.

For time-domain problems, as the discontinuity moves, one has to efficiently update these correction terms in each time-step. 
Substituting Eq.~\eqref{sjofx2b} into \eqref{fnderivativeapprox2} yields the suggestive expression
\begin{equation}
\label{fnderivativeapprox3}
{f^{(n)}}_i \simeq \sum\limits_{j = 0}^N D_{ij}^{(n)}(f_j+
 g_j \theta_{i}-g_j \theta_j),
\end{equation}
where $f_i=f(x_i), {f^{(n)}}_i={f^{(n)}}({x_i})$, 
$\theta_i=\theta (x_{i} - \xi )$ and $g_j=g(x_j-\xi)$ are vectors formed from the values of the respective functions on the set of grid-points $\{x_i\}$.
On modern CPUs and GPUs, Eq.~\eqref{fnderivativeapprox3} can be evaluated efficiently via (inner) matrix-vector and (element-wise) vector-vector multiplication.
For instance, in Wolfram Language, the simple command
\begin{verbatim}
fn = Dn.f + (Dn.g)*h - Dn.(g*h)
\end{verbatim}
(with \verb|f| and \verb|fn|  denoting $f_i$ and ${f^{(n)}}_i$,
\verb|Dn| denoting the $n^{\rm{th}}$ order differentiation matrix given by Eq.~\eqref{fnderivativeapprox1}, and \verb|h|, \verb|g| denoting the vectors  $\theta_i$, $g_i$ respectively)
%
uses the \textsc{Intel Math Kernel Library (Intel MKL)} to automatically perform the linear algebra operations in Eq.~(\ref{fnderivativeapprox3}) in parallel, accross all available cores. 
With sufficient optimization and hyperthreading enabled, our implementation saturates all 28 cores of a Intel Xeon  8173M CPU ($\sim 100\%$ core utilization). 
 A similar code, that loads the relevant matrices onto GPU memory and replaces the \verb .  or 
\verb Dot[]  command with \verb CUDADot[] , uses the \textsc{CUDA BLAS (cuBLAS)} Library to accelerate these linear algebra operations on Nvidia GPUs.
In our experiments, for sufficiently large matrices, Nvidia Volta architecture GPUs with double precision support, such as a Titan V or  Quadro GV100, provided a $\sim 5 \times$ speed-up over a CPU of similar Wattage, such as a 28-core Intel Xeon 8173M. OpenCL acceleration on CPUs or GPUs is also possible.

As a finite-difference example, let us consider a function with jumps $J_0, J_1$ and $J_2$ in its $0^{\rm{th}}$, $1^{\rm{st}}$ and $2^{\rm{nd}}$ derivative at the  $J^{\rm{th}}$ node, $x_J=\xi$. Then, for $3$-point  finite-difference stencils, Eq.~\eqref{fnderivativeapprox2}
gives
the left and right derivatives:
\begin{eqnarray} \label{fnderivativeleftright}
f'_{J^{\pm}}= f'_J &-& \frac{1}{2}{J_0}\left( {\frac{1}{{{x_J} - {x_{J - 1}}}} + \frac{1}{{{x_{J + 1}} - {x_J}}} - \frac{2}{{{x_{J + 1}} - {x_{J - 1}}}}} \right)\\
&+&J_1\frac{{{x_{J \pm 1}} - {x_J}}}{{{x_{J + 1}} - {x_{J - 1}}}}
- \frac{1}{2}{J_2}
\frac{{({x_{J + 1}} - {x_J})({x_J} - {x_{J - 1}})}}{{{x_{J + 1}} - {x_{J - 1}}}},
\nonumber
\end{eqnarray}
with $f'_J$  given by Eq.~\eqref{smoothFDf} for $i= J$.
Notice that the above formula satisfies 
the jump condition
 $f'_{J^{+}}-f'_{J^{-}}=J_1$ 
 as it ought to.

 Higher order finite-difference or pseudo-spectral formulas can be constructed in a similar fashion.
With the nodes given by Eq.~\eqref{ChebyshevNodes}, 
the formula \eqref{fnderivativeapprox2}
yields a discontinuous collocation method of $M^{\rm {th}}$-order convergence  if  derivatives of $f$ higher than $M$ have finite jumps, or of infinite-order convergence  if derivatives higher than $M$ are continuous.
This will be demostrated explicitly below.

\subsection{Discontinuous spatial integration}

Numerical integration formulas for  nonsmooth integrands can be obtained, as in \S2.3,  by integrating the interpolating formula \eqref{PolyNLagrange2}. Then, the standard quadrature formula  
 \eqref{Integralp1} generalizes to
\begin{equation} \label{Integralp2}
\int_{{a}}^{{b}} {f(x)dx}\simeq  \int_{{a}}^{{b}} {p(x)dx}  = \sum\limits_{j = 0}^N {{[w_{j}}{f_j}+q_j(\xi)]},
\end{equation}
where the weights
$w_j$ are given by Eq. \eqref{IntegralMatrixp1}, and the corrections
\begin{equation} \label{IntegralMatrixp2}
{q_{j}(\xi)}=\int_{{a}}^{{b}} {\Delta(x_j-\xi;x-\xi)\pi_j(x)dx}
\end{equation}
account for the discontinuity and vanish when $J_m=0$.
Equation~\eqref{Integralp2} can be used to obtain discontinuous generalizations of standard composite or Gaussian quadrature formulas using, for example, the nodes \eqref{EquidistantNodes} or \eqref{ChebyshevNodes} respectively.  
This formula can be useful for spatial integration of discontinuous or non-smooth functions, that arise, for instance, when one wishes to compute field energy or angular momentum integrals in the presence of a moving particle over a spatial slice. For the purpose of time integration, however, we will use a 2-point time-symmetric method, also based on the method of undetermined coefficients, detailed below.

\begin{figure}[ht] 
\includegraphics[width=\columnwidth]{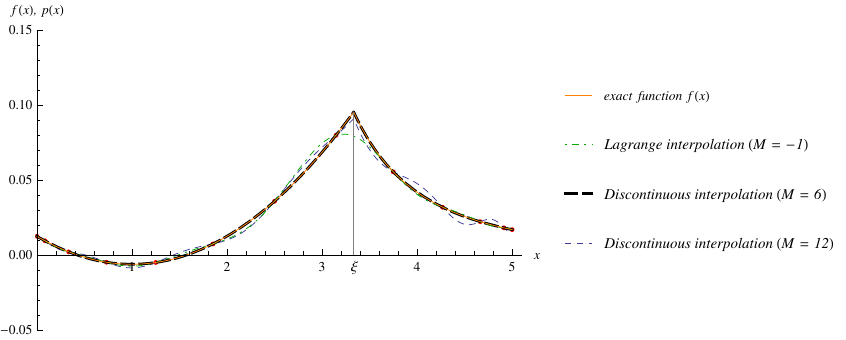}
\includegraphics[width=\columnwidth]{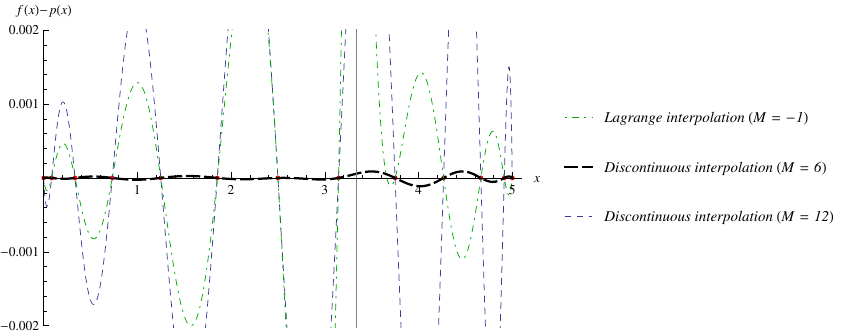}
\caption{Same as Fig. 1 but using Chebyshev-Gauss-Lobatto nodes given by Eq.~\eqref{ChebyshevNodes},
which minimize the maximum error  by distributing oscillations more evenly. 
Eq.~\eqref{PolyNLagrange2}    faithfully reproduces the function near and away from the discontinuity, especially when an intermediate number  ($M\approx N/2$) of derivative jumps are enforced.} 
\label{diff2} 
\end{figure}  

\begin{figure}[ht] 
\includegraphics[width=\columnwidth]{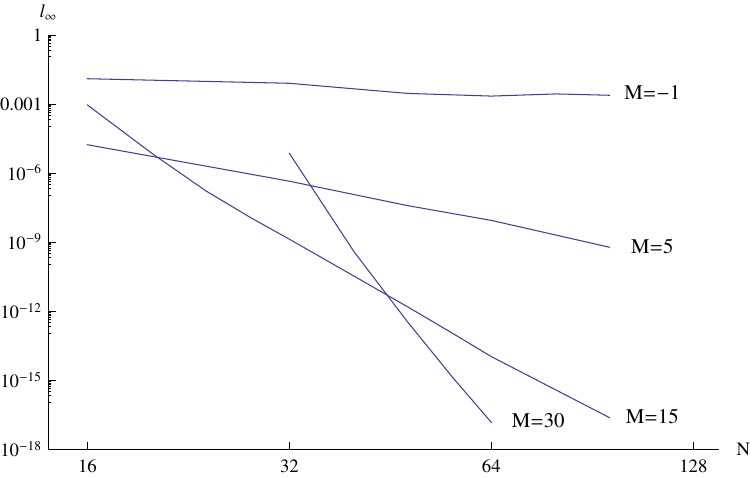}
\caption{$l_\infty$ error norm of interpolating the function ${f}(x) = {P_2}(\xi ){Q_2}(x)\theta (x - \xi ) + {P_2}(x){Q_2}(\xi )\theta (\xi - x)$  at $N+1$ Chebyshev-Gauss-Lobatto  nodes using the Lagrange formula~\eqref{PolyNLagrange1} or its discontinuous generalization \eqref{PolyNLagrange2} with jump conditions  on the first $M=5, 15$ or $30$ derivatives enforced. Having dropped the (nonzero in this case) jumps in derivatives  higher   than $M$,  we observe convergence of  order slightly above $M$. For functions with continuity in derivatives higher than $M$,    infinite-order (exponential) convergence  is attained.} 
\label{diff3} 
\end{figure}  

\section{Application to PDEs with distributional sources} \label{sec:sec4}
A prime application of the discontinuous collocation methods outlined above is the numerical evolution of partial differential equations whose solutions are known to be non-analytic at a single point. These problems typically arise in self-force computations and black hole perturbation theory. The discontinuous behavior can be captured by splitting the computational domain into two regions with the particle's worldline as the time-dependent boundary. 
The relevant PDEs are then solved in time-dependent (Lagrangian) coordinates that track the particle motion.
We will now demonstrate that discontinuous collocation methods allow these problems to be solved without domain decomposition, in fixed (Eulerian) coordinates, with the contribution of the moving particle captured by known time-dependent terms linear in the (known) jumps. This approach {\it{(i)}} reduces code complexity, {\it{(ii)}} provides users with the flexibility of selecting their own spacetime coordinates (such as hyperboloidal slicing) that accommodate their specific computational needs,
{\it{(iii)}}
avoids the need of expensive re-gridding in each time-step, 
{\it{(iv)}}
allows conversion of finite-difference to pseudo-spectral discretizations, or extension to higher order, by simply changing one code parameter, rather than completely rewriting the code, and
{\it{(v)}}
allows for straightforward parallelization in multi-core CPU and GPGPU architectures.

\subsection{Discontinuous time integration}
The methods outlined above enable us to easily account for the non-analytic behavior of the target function at a single point along the spatial axis. However, the function is also non-analytic when the point is approached along the temporal axis. This indicates that standard time-steppers (e.g. Runge-Kutta methods) cannot be readily applied to this problem, as they assume the target function is smooth.
Our numerical (method-of-line) experiments indeed showed that using ordinary time-steppers with the (spatially) discontinuous collocation method described above works accurately if the particle is static, but precision is lost if the particle is moving.
We will show here that the paradigm of using undetermined coefficients to accommodate known jump and collocation conditions can also be used to develop discontinuous time integration schemes. 

For concreteness, we consider the differential equation
\begin{equation}
    \frac{dy}{dt} = f(y,t)
\end{equation}
and approximate its solution on a small time interval $[t_1,t_2]$. This is equivalent to writing
\begin{equation}
    y(t_2) - y(t_1) = \int_{t_1}^{t_2} f(y,t) dt,
\end{equation}
so the problem is now to approximate the above integral. One approach is to construct a Lagrange polynomial approximant to $f$ between $t_1$ and $t_2$. Since there exist  two collocation points, the interpolating polynomial is of degree one. Applying Eq. \eqref{PolyNLagrange1}, we obtain the linear interpolating function
\begin{equation}
    f(y,t) \simeq p(y,t) = f(y_1, t_1) ~ \frac{t_2 - t}{t_2 - t_1} + f(y_2,t_2) ~ \frac{t - t_1}{t_2 - t_1}.
\end{equation}
Inserting this expression into the above integral, we obtain the trapezium rule
\begin{equation}
    y_2 - y_1 \simeq \frac{\Delta t}{2} \big( f_1 + f_2 \big ).
\end{equation}
In a method-of-lines context, when combined with second-order (spatial) finite differencing, this amounts to the Crank-Nicolson method of solving time-dependent PDEs. We have thus constructed an implicit integration scheme \footnote{There are several theoretical advantages to implicit time-symmetric methods, such as unconditional stability and conservation of energy \cite{2019arXiv190109967M}.}
by applying smooth interpolation methods to the differential equation.

If the function $y$ (and by extension $f$) or its derivatives have jump discontinuities at a time $t=\tau \in [t_1,t_2]$, then the standard Lagrange interpolating polynomial approximation fails to converge pointwise near $\tau$. It is, however, possible to proceed by the method used in Section \ref{InterpolationSection}. The interval is split into two parts and $f$ is approximated by $p_<$ and $p_>$ in $[t_1,\tau]$ and $[\tau,t_2]$ respectively:
\begin{align}
    p_<(t) &= A + B t\\
    p_>(t) &= C + D t.
\end{align}
The above polynomials together approximate the discontinuous function $f$ over the interval $[t_1,t_2]$. They can be regarded together as a single piecewise polynomial
\begin{equation}
p(t)=\theta(t-\tau) p_> (t)+\theta(\tau-t) p_< (t). 
\end{equation}
There are four unknown coefficients along with four conditions: the collocation conditions $p(t_1) = f_1$, $p(t_2) = f_2$
and the jump conditions
$p(\tau^+) - p(\tau^-) = J_0$, $p'(\tau^+) - p'(\tau^-) = J_1$.
The four conditions amount to a linear system for the four polynomial coefficients. Solving the system yields:
\begin{equation}\label{pGreat}
    p_> (t) = f_1 ~ \frac{t_2 - t}{t_2 - t_1} + f_2 ~ \frac{t - t_1}{t_2 - t_1} + \frac{t_2 - t}{t_2 - t_1}( J_0 - J_1 (\tau - t_1) )
\end{equation}
\begin{equation}\label{pLess}
    p_< (t) = f_1 ~ \frac{t_2 - t}{t_2 - t_1} + f_2 ~ \frac{t - t_1}{t_2 - t_1} + \frac{t - t_1}{t_2 - t_1}( - J_0 - J_1 (t_2 - \tau) ).
\end{equation}
Integrating this piecewise polynomial over the time interval $[t_1,t_2]$ yields the discontinuous generalization of the trapezium rule:
\begin{equation}\label{eq:discTrap}
    y_2-y_1 \simeq  \frac{\Delta t}{2} ( f_1 + f_2) + J_0 ~ \frac{\Delta t - 2 \Delta \tau}{2} + J_1 ~ \frac{\Delta \tau}{2} ( \Delta \tau - \Delta t ),
\end{equation}
where $\Delta \tau = \tau - t_1$. The standard trapezoidal rule is obviously recovered when $J_0 = J_1 = 0$. We have thus constructed an implicit integration scheme that allows the target function $y$ to have non-analytic behavior. Although this analysis has  yielded a method second-order  in time, it  is fairly straightforward to extend to higher order,  by using higher order piecewise (Lagrange or Hermite) interpolating polynomials in each time interval  \cite{2019arXiv190109967M} or composing this method to create a higher order multi-stage scheme \cite{2010Lubina}. A fourth-order discontinuous Hermite rule for time integration accross the particle worldline will be  demonstrated  in a companion paper. We now apply the discontinuous trapezium rule  to a simplified version of the time-dependent self-force problem.

\subsection{Distributionally-forced wave equation}\label{sec:deltaWave}
We will study a simplified version of Eq.~\eqref{TeukSchw} presented in \cite{Field2009},
\begin{equation}\label{eq:Field1}
    - \partial^2_t \Phi + \partial^2_x \Phi = \cos t ~ \delta( x - v t ).
\end{equation}
One can perform a first order reduction in time by defining $\Pi = \partial_t \Phi$:
\begin{align}
    \label{reducedPhi}
    \partial_t \Phi &= \Pi\\
    \label{reducedPi}
    \partial_t \Pi & = \partial^2_x \Phi - \cos{t} ~ \delta( x - v t ).
\end{align}
We handle the discontinuous behavior in $\Phi$ (and by extension $\Pi$) by approximating it with a piecewise polynomial
\begin{equation}\label{discreteDerivative}
\Phi(t,x) \simeq \sum\limits_{j = 0}^N {{[\Phi_j(t)+\Delta(x_j-\xi(t);x-\xi(t))]}{\pi _j}(x)}.
\end{equation}
Recall that, by construction, the functions $\Delta(x_j-\xi;x-\xi)$ are piecewise constant with respect to $x$, so they are not differentiated when approximating the spatial derivatives:
\begin{equation}
    \partial_x^n \Phi (t, x) |_{x=x_i} \simeq \sum\limits_{j = 0}^N D_{ij}^{(n)}[\Phi_j+\Delta(x_j-\xi;x_i-\xi)].
\end{equation}
with the differentiation matrices given by Eq.~(\ref{derivativematricespi}). There is no summation over $i$, so, the rightmost term in the above expression is a vector of length $N+1$, denoted by
\begin{equation}\label{eq:effectiveSource}
    s^{(n)}_i (t) = \sum\limits_{j = 0}^N D_{ij}^{(n)}\Delta(x_j-\xi(t);x_i-\xi(t)).
\end{equation}
The details of calculating the jump terms necessary for the construction of $s$ are given in \ref{appendix:jumps}. The full discretized system accounting for discontinuities and boundary conditions is derived in \ref{appendix:BCs} and \ref{appendix:MOL}. The discretized version of Eqs. \eqref{reducedPhi} and \eqref{reducedPi} is given in Eqs. \eqref{discretePhi} and \eqref{discretePi}.
\begin{figure}[t]
    \centering
    \includegraphics[scale=0.9]{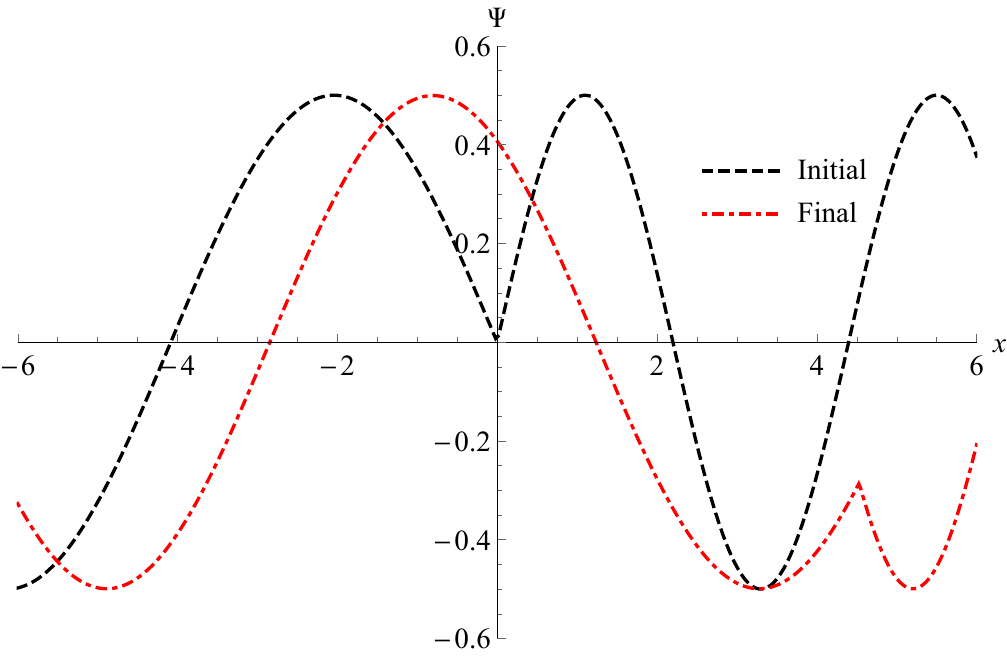}
    \caption{Numerical solution to Eq. \eqref{eq:Field1} using the discontinuous collocation methods of §3 for spatial discretization and the discontinuous trapezoidal rule of §4 for temporal integration. The Initial function was selected to be the exact solution provided in \cite{Field2009} at $t=0$; this ensures the numerical solution can be directly compared to an exact solution when studying convergence. We evolved the reduced system from $t=0$ to $t=9.5$ taking $\Delta t = 0.01$.}
    \label{fig:Field_Sol}
\end{figure}

\begin{figure}[t]
    \centering
    \includegraphics[scale=0.7]{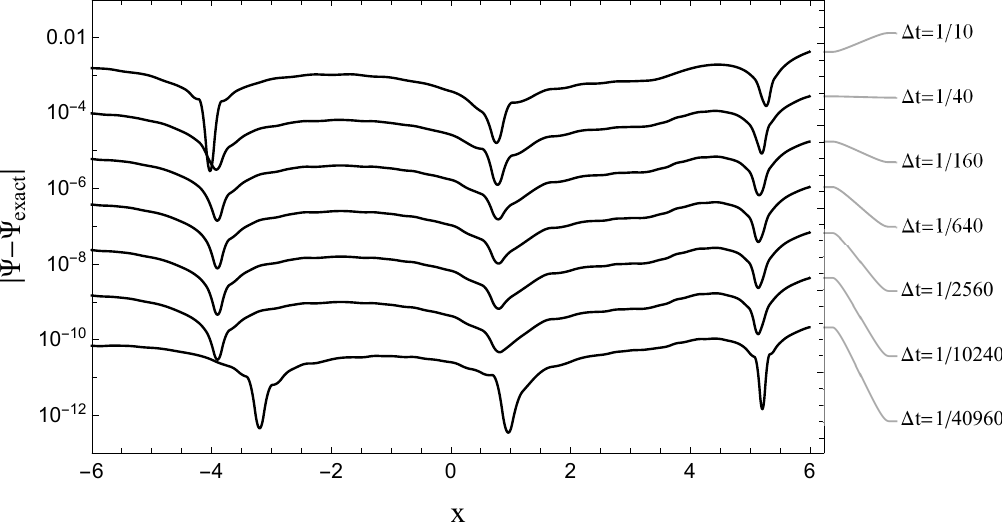}\\
    \includegraphics[scale=0.7]{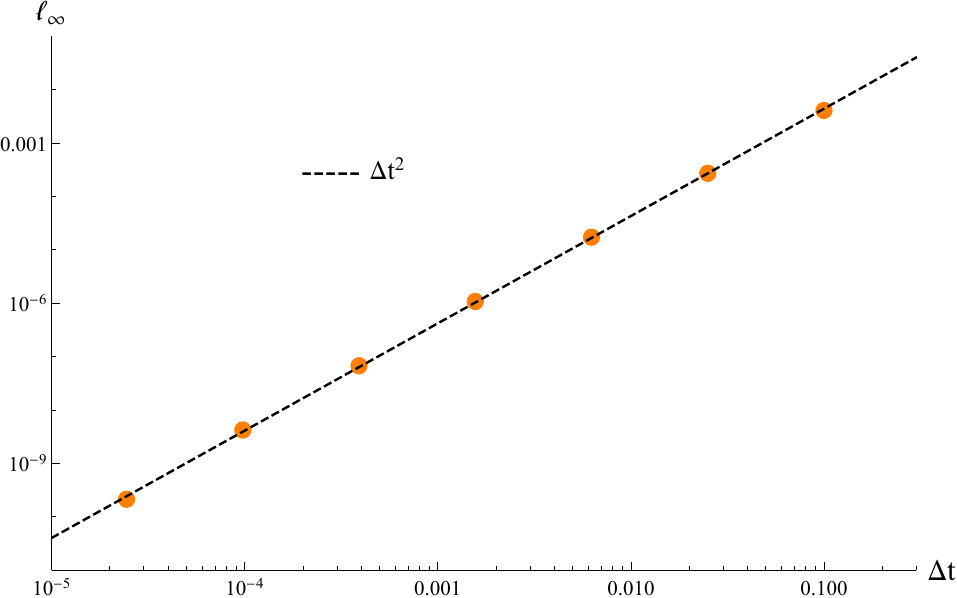}
    \caption{The top figure shows the difference between the numerical solution and the exact solution to Eq. \eqref{eq:Field1} for decreasing time steps. The bottom figure demonstrates that the $l_\infty$ error norm scales as $\Delta t^2$, as expected of our method.}
    \label{fig:TrapError}
\end{figure}

We apply the discontinuous version of the trapezoidal rule derived in the previous section to numerically evolve this problem in time. We take $v=0.3$ in Eq. \eqref{reducedPi} and consider a Chebyshev pseudospectral grid with $a = -6$ and $b = 6$ and $110$ grid points. The authors of \cite{Field2009} provided an exact solution to this problem, which our numerical scheme reproduces as shown in Figure \ref{fig:Field_Sol}. We found that including $M\simeq30$ spatial jumps was sufficient. Moreover, we observe second order convergence in time as is illustrated in Figure \ref{fig:TrapError}.

As demonstrated in Ref.~\cite{2019arXiv190109967M}, the advantages of the implicit time-symmetric integration method used here (as opposed, for instance, to explicit Runge-Kutta methods) are that {\it{(i)}} energy and symplectic struncture are conserved exactly for linear ODEs and PDEs (and on average for nonlinear ODEs and PDEs), {\it{(ii)}} symmetric schemes are unconditionally stable, that is, not subject to a Courant–Friedrichs–Lewy condition, and {\it{(iii)}} the implicitness of the scheme does not amount to extra computational cost for the linearized PDEs arising in black hole pertubation theory, as linearized differential operations can be descretized, inverted and stored in memory. These are major advantages that can be exploited for long-time numerical evolutions, necessary for modelling EMRIs over the 2-5 years they are expected to spend in LISA band. In future work it will be shown that a discontinuous Hermite or Lotkin integration method can be used to obtain fourth  or sixth order convergence in time.

\section{Summary and applications}

We have shown that, if the  derivative jumps across a discontinuity are known, the simple modification \eqref{PolyNLagrange2} to the Lagrange formula \eqref{PolyNLagrange1}  can  fully  account for the presence of a discontinuity and completely avoid Gibbs phenomena.
This can be exploited in finite-difference or pseudospectral scheme to recover their original accuracy. The main advantage offered by the present method is that it is efficient and simple to implement, even when the discontinuity  moves freely inside the domain, as no requirements for domain decomposition, regridding, or reconstruction arise.  

Accurate numerical differentiation of piecewise smooth functions is crucial to many problems of physical interest.  A common scenario consists of evolving  a time-dependent partial differential equation via the method of lines \cite{knapp2008method,Schiesser1991}. Spatial   operators are then computed via  finite-difference or pseudospectral methods described above, and time integration  is performed using  standard numerical methods  for coupled ordinary differential equations, such as Runge-Kutta. The standard methods of spatial discretization in \S2 cease to apply when discontinuities are present; this is when Eq.~\eqref{fnderivativeapprox2} becomes useful. For a moving discontinuity, the weights \eqref{sjofx2b} must by evaluated at each time step; if the jumps $\{J_m\}$ at $\xi$ are known, this evaluation is not  expensive computationally as it requires a few subtractions and multiplications but no  extra function evaluations.
The computational cost and accuracy of obtaining the jumps and their location depends on the problem at hand; in certain problems  this  information is readily available at no cost.
In the case of self-force applications such as that presented above, the jumps \eqref{fmJconditions2}
occur at the known particle location and can be computed analytically from the radial part of the field equations.
The present method is ideally suited to such problems where information on the jumps is available.

Lipman and Levin \cite{Lipman2009} have obtained an   interpolation formula for piecewise smooth functions written in a monomial basis, which can be considered equivalent to 
Eq.~(\ref{PolyNLagrange2}) above written in a Lagrange basis. Introducing a Lagrange basis is essential, 
  as it allowed us to go beyond interpolation and construct differentiation matrices for discontinuous collocation (finite-difference or pseudo-spectral) methods  demonstrated
above. Lipman and Levin moreover devised a moving least-squares technique for detecting the location and magnitude of  an unknown discontinuity or singularity, with accuracy of higher order  than subcell resolution methods \cite{Harten1989}. This least-squares technique allows use of the methods presented here to problems with
time-dependent discontinuities that are not known a priori. 

 For  pseudospectral discretizations,   a mapping devised by Kosloff and Tal-Ezer 
 to improve  stability 
of the Chebyshev collocation scheme 
\cite{Kosloff1993} and  the accuracy of higher derivatives \cite{Costa2000,Don1997}, can  also be employed
in conjunction with  Eq.~(\ref{fnderivativeapprox2}) in  method-of-lines implementations to evolve  discontinuous solutions.

In problems where the grid values of the function derivatives    are known during evolution (in addition to their jumps at the discontinuity),
Eqs.~(\ref{PolyNConds2}) can be complemented by conditions on the derivatives at the same nodes. Together  with Eq.~(\ref{fmJ0conditions2}), this will yield generalizations of Hermite \cite{hermite1878}
and Birkhoff \cite{Birkhoff1906} interpolation to discontinuous problems. This can be very useful for increasing the  order of our discontinuous time integrator, as will be shown in a companion paper.

As a final remark, we note that extension of the present technique to higher spatial dimension is technically possible, by replacing Eq.~\eqref{gjofx2} with a multi-dimensional Taylor expansion. 
Although such an expansion would in principle work,
the number of mixed derivative jump terms would increase very rapidly with the order $M$ of the method, which making it prohibitive to go to very high order.
However, for separable problems, mixed derivative jump terms are absent. For instance, a (spin-weighted) spherical harmonic expansion can be used to eliminate the angular terms in the Teukolsky \cite{Csukas2019kcb,Dolan2013,Racz2011} or Regge-Wheeler-Zerilli \cite{Field2009,Field:2010xn} equations and reduce them to 1+1 equations with radial and time dependence only. We thus believe that our method is well-suited to near-separable problems on tensor-product grids \cite{2020Wolfram,knapp2008method}.

%

Some of these ideas are currently being implemented in a numerical scheme where the metric perturbation is reconstructed in the time domain, by directly evolving a sourced Teukolsky equation satisfied by the Hertz potential \cite{Giudice2015}.

\section*{Acknowledgments}
We thank Sarp Akcay, Bernd Br\"ugmann, Priscilla Ca\~nizares, Peter Diener, Sam Dolan, Scott Field, Enno Harms, Ian Hawke, David Hilditch, Jonathan Thornburg, Maarten van de Meent, Leo Stein and K\=oji Ury\=u for constructive discussions and comments. The research leading to these
results received funding from the European Research Council under the European Union's
Seventh Framework Programme (FP7/2007-2013)/ERC grant agreement no. 304978
and the European Union's Horizon 2020 research and innovation programme under the Marie Sk\l odowska-Curie grant agreement no. 753115.
This work was further supported by STFC grant PP/E001025/1 and DFG grant SFB/Transregio 7 ``Gravitational Wave Astronomy''.

\appendix
\section{Calculating Jumps}\label{appendix:jumps}
The interpolation function of Eq. \eqref{PolyNLagrange2} requires knowledge of the discontinuities in the first $M$ spatial derivatives of the target function $f$. However, in self-force calculations (or, more generally, fields sourced by a point charge), the underlying equation only specifies the discontinuities in the zeroth and first spatial derivatives. We will now demonstrate that jump in the spatial derivative of any order is determined by the first two jumps.

As a simple example, consider a one-dimensional wave equation sourced by a singularity moving with trajectory $\xi(t)$:
\begin{equation}
    - \frac{\partial^2 \Phi}{\partial t^2} + \frac{\partial^2 \Phi}{\partial x^2} = F(t) \delta( x - \xi(t) ) + G(t) \delta'(x - \xi(t))
\end{equation}
The general solution $\Phi(x,t)$ may be decomposed into a smooth part representing the homogeneous solution and a non-smooth part capturing the jumps. We will consider only the non-smooth part $\Phi_J$.

We define ``unit jump functions'' by
\begin{equation}
    \Psi_n (x; \xi) = \frac{1}{2} \mathrm{sgn}(x - \xi) \frac{(x - \xi)^n}{n!}
\end{equation}
where $\mathrm{sgn}$ is the signum function. These functions have the property that all spatial derivatives have matching left and right limits at $x=\xi$ except the $n^\mathrm{th}$ derivative, which has a jump of one. That is,
\begin{equation}
    \frac{\partial^m \Psi_n}{\partial x^m} \bigg\rvert_{\xi^{+}} - ~ \frac{\partial^m \Psi_n}{\partial x^m} \bigg\rvert_{\xi^{-}} = ~ \delta_{m,n}
\end{equation}
So, if the discontinuity in the $n^\mathrm{th}$ derivative of $\Phi$ is given by $J_n(t)$, we may express $\Phi_J$ as a linear combination of these jumps
\begin{equation}
    \Phi_J(x,t) = \sum_{i=0}^\infty J_i(t) \Psi_i (x;\xi(t))
\end{equation}

This ansatz is substituted into the original equation and terms of the same order are grouped. The following relations emerge:
\begin{align}
    J_0 =& \gamma^2 G(t)\\
    J_1 = \gamma^2( F(t) &- 2 \Dot{\xi} \Dot{J}_0 - \Ddot{\xi} J_0)\\
    J_{n+2} = \gamma^2( \Ddot{J_n} -& 2 \Dot{\xi} \Dot{J}_{n+1} - \Ddot{\xi} J_{n+1})
\end{align}
where
\begin{equation}
    \gamma^2 = \frac{1}{1 - \Dot{\xi}^2}
\end{equation}
Thus, the original equation directly determines $J_0$ and $J_1$, and this specifies all higher order jumps in the spatial derivative.

\section{Radiation Boundary Conditions}\label{appendix:BCs}
The spatial domain of Eq. \eqref{eq:Field1} and its solution presented in \cite{Field2009} is the entire real line. Care must be taken when numerically solving this problem on a finite computational domain. One method that has gained popularity in self-force applications is to use a spatial compactification and hyperboloidal time slicing; this way, the behavior at ``physical infinity'' is represented at a finite coordinate in finite time \cite{2011JCoPh.230.2286Z}.

In this work, we chose to follow \cite{Field2009} and apply radiation (or Sommerfeld) boundary conditions to our domain. The wave equation has two characteristic speeds; in one dimension, they correspond to left- and right-propagating waves. However, only one of these speeds is physical at each of the boundaries. For example, a right-propagating wave at $x=a$ corresponds to radiation from outside sources entering the domain, which is unphysical to our model. This problem is avoided if we require solutions to be entirely left-propagating at $x=a$ (or purely ``outgoing''). That is, we require
\begin{equation}\label{eq:leftWave}
    \Pi(x=a) = \partial_x \Phi(x=a)
\end{equation}
The corresponding condition on the right boundary is
\begin{equation}\label{eq:rightWave}
    \Pi(x=b) =- \partial_x \Phi(x=b)
\end{equation}

It is possible to impose these conditions by expressing them as algebraic constraints on the discrete grid. Recalling $x_0=a$ and $x_N=b$, we have
\begin{equation}
    \Pi_0 = \sum_{j=0}^N D^{(1)}_{0 j} \Phi_j + s_0^{(1)}(t)
\end{equation}
\begin{equation}
    \Pi_N = -\sum_{j=0}^N D^{(1)}_{N j} \Phi_j - s_N^{(1)}(t)
\end{equation}
The function $s_i^{(1)}(t)$ is the same as defined in Eq. \eqref{eq:effectiveSource}. We proceed by solving for $\Phi_0$ and $\Phi_N$ and eliminating them from the discretized problem. However, additional constraints are needed, as $\Pi_0$ and $\Pi_N$ must also be eliminated. This is obtained by taking spatial derivatives of Eqs. \eqref{eq:leftWave} and \eqref{eq:rightWave}, so we also have
\begin{equation}
    \sum_{j=0}^N D^{(1)}_{0 j} \Pi_j + r_0^{(1)}(t) = \sum_{j=0}^N D^{(2)}_{0 j} \Phi_j + s_0^{(2)}(t)
\end{equation}
\begin{equation}
    \sum_{j=0}^N D^{(1)}_{N j} \Pi_j + r_N^{(1)}(t) = -\sum_{j=0}^N D^{(2)}_{N j} \Phi_j - s_N^{(2)}(t)
\end{equation}
$r_i^{(n)}(t)$ is defined in the same manner as Eq. \eqref{eq:effectiveSource} but is evaluated using the jumps in $\Pi$ as opposed to the jumps in $\Phi$. We may now solve for $\Phi_0$ and $\Phi_N$:
\begin{multline}
    \Phi_0 = \frac{1}{\Omega} \sum_{j=1}^{N-1} \bigg( (A_0 D_{0j}^{(1)} + B_0 D_{Nj}^{(1)} + C_0 D^{(2)}_{0j} + E_0 D^{(2)}_{Nj}) \Phi_j + (-C_0 D_{0j}^{(1)} + E_0 D_{Nj}^{(1)}) \Pi_j \bigg)\\
    + A_0 s_{0}^{(1)}(t) + B_0 s_{N}^{(1)}(t) + C_0 s^{(2)}_{0}(t) + E_0 s^{(2)}_{N}(t) - C_0 r_{0}^{(1)}(t) + E_0 r_{N}^{(1)}(t)
\end{multline}
\begin{multline}
    \Phi_N = \frac{1}{\Omega} \sum_{j=1}^{N-1} \bigg( (A_N D_{0j}^{(1)} + B_N D_{Nj}^{(1)} + C_N D^{(2)}_{0j} + E_N D^{(2)}_{Nj}) \Phi_j + (-C_N D_{0j}^{(1)} + E_N D_{Nj}^{(1)}) \Pi_j \bigg)\\
    + A_N s_{0}^{(1)}(t) + B_N s_{N}^{(1)}(t) + C_N s^{(2)}_{0}(t) + E_N s^{(2)}_{N}(t) - C_N r_{0}^{(1)}(t) + E_N r_{N}^{(1)}(t)
\end{multline}
where
\begin{align*}
    \Omega &= (D_{0N}^{(1)})^2 (D_{N0}^{(1)})^2 - (D_{NN}^{(1)})^2 D^{(2)}_{00} - D_{00}^{(1)} D_{N0}^{(1)} D_{0N}^{(2)} + D_{N0}^{(1)} D_{NN}^{(1)} D_{0N}^{(2)} - D_{0N}^{(2)} D_{N0}^{(2)} \\
    &+ (D_{00}^{(1)})^2 ( (D_{NN}^{(1)})^2 - D_{NN}^{(2)}) + D_{00}^{(2)} D_{NN}^{(2)} + D_{0N}^{(1)}(-D_{NN}^{(1)} D_{N0}^{(2)} \\
    &+ D_{00}^{(1)} ( D_{N0}^{(2)} - 2 (D_{N0}^{(1)} D_{NN}^{(1)}) + D_{N0}^{(1)} (D_{00}^{(2)} + D_{NN}^{(2)})))\\\\
    A_0 &= D_{0N}^{(1)} D_{N0}^{(1)} D_{NN}^{(1)} - D_{00}^{(1)} (D_{NN}^{(1)})^2 + D_{N0}^{(1)} D_{0N}^{(2)} + D_{00}^{(1)} D_{NN}^{(2)}\\
    B_0 &= D_{0N}^{(1)} D_{00}^{(1)} D_{NN}^{(1)} - (D_{0N}^{(1)})^2 D_{N0}^{(1)} - D_{NN}^{(1)} D_{0N}^{(2)} - D_{0N}^{(1)} D_{NN}^{(2)}\\
    C_0 &= (D_{NN}^{(1)})^2 - D_{0N}^{(1)} D_{N0}^{(1)} - D_{NN}^{(2)}\\
    E_0 &= D_{0N}^{(1)} D_{NN}^{(1)} - D_{00}^{(1)} D_{0N}^{(1)} + D_{0N}^{(2)}\\\\
    A_N &= D_{N0}^{(1)} D_{00}^{(1)} D_{NN} - (D_{0N}^{(1)})^2 D_{N0}^{(1)} - D_{N0}^{(1)} L_{00}^{(2)} - D_{00}^{(1)} D_{N0}^{(2)}\\
    B_N &= D_{00}^{(1)} D_{0N}^{(1)} D_{N0}^{(1)} - (D_{00}^{(1)})^2 D_{NN}^{(1)} + D_{NN}^{(1)} D_{00}^{(2)} + D_{0N}^{(1)} D_{N0}^{(2)}\\
    C_N &= D_{00}^{(1)} D_{N0}^{(1)} - D_{N0}^{(1)} D_{NN}^{(1)} + D_{N0}^{(2)}\\
    E_N &= (D_{00}^{(1)})^2 - D_{0N}^{(1)} D_{N0}^{(1)} - D_{00}^{(2)}
\end{align*}
These expressions for $\Phi_0$ and $\Phi_N$ may now be substituted into Eq. \eqref{discreteDerivative} to obtain
\begin{equation}
    \sum_{j=0}^N D^{(2)}_{ij} \Phi_j + s_i^{(2)}(t) \rightarrow \sum_{j=1}^{N-1} \Bar{D}_{ij}^{(2)} \Phi_j + \sum_{j=1}^{N-1} M_{ij} \Pi_j + \Bar{s}_i(t)
\end{equation}
where $i$ is restricted to $1,2,\dots,N-1$ and
\begin{align*}
\Bar{D}_{ij}^{(2)} &= D_{ij}^{(2)} + \frac{1}{\Omega} D_{i0}^{(2)} (A_0 D_{0j}^{(1)} + B_0 D_{Nj}^{(1)} + C_0 D^{(2)}_{0j} + E_0 D^{(2)}_{Nj} ) \\
   &  + \frac{1}{\Omega} D_{iN}^{(2)}(A_N D_{0j}^{(1)} + B_N D_{Nj}^{(1)} + C_N D^{(2)}_{0j} + E_N D^{(2)}_{Nj})\\
M_{ij}&=\frac{1}{\Omega} D_{i0}^{(2)}(-C_0 D_{0j}^{(1)} + E_0 D_{Nj}^{(1)} ) + \frac{1}{\Omega} D_{iN}^{(2)}(-C_N D_{0j}^{(1)} + E_N D_{Nj}^{(1)} )   \\
\Bar{s}_i(t) &= s^{(2)}_i(t) + D_{i0}^{(2)} \Big( A_0 s_{0}^{(1)}(t) + B_0 s_{N}^{(1)}(t) + C_0 s^{(2)}_{0}(t) + E_0 s^{(2)}_{N}(t) - C_0 r_{0}^{(1)}(t) + E_0 r_{N}^{(1)}(t) \Big)\\
   & + \Bar{D}_{iN}^{(2)} \Big( A_N s_{0}^{(1)}(t) + B_N s_{N}^{(1)}(t) + C_N s^{(2)}_{0}(t) + E_N s^{(2)}_{N}(t) - C_N r_{0}^{(1)}(t) + E_N r_{N}^{(1)}(t) \Big)
\end{align*}
%
We see that the boundary conditions  \eqref{eq:leftWave} and \eqref{eq:rightWave} are accommodated by removing the first and last points from the spatial grid, modifying the second derivative operator acting on the field $\Phi$, and introducing coupling to the momentum $\Pi$ and a new effective source term. 

\section{Method-of-Lines Discretization}\label{appendix:MOL}
We discuss how the reduced system Eqs. \eqref{reducedPhi} and \eqref{reducedPi} may be approximately integrated in time using the discontinuous trapezium rule Eq. \eqref{eq:discTrap}. We begin by discretizing Eqs. \eqref{reducedPhi} and \eqref{reducedPi} in space in the manner outlined in \ref{appendix:BCs} accounting for boundary conditions and discontinuities. Then, our partial differential equations amount to a systems of ordinary differential equations coupled via differentiation matrices:
\begin{align}
    \label{spacePhi}
    \frac{d \Phi_i}{dt} &= \Pi_i\\
    \label{spacePi}
    \frac{d \Pi_i}{dt}  &= \Bar{D}_{ij}^{(2)} \Phi_j + M_{ij} \Pi_j + \Bar{s}_i(t) - \cos{t}~\delta(x_i - v t)
\end{align}
We use the Einstein summation convention in this section, implying summation over the index $j$ (but not $i$), which amounts to matrix-vector multiplication.
These equations are then approximately integrated in time using the discontinuous trapezium rule \eqref{eq:discTrap}. Applied to Eq. \eqref{spacePhi}, we have
\begin{equation}\label{discretePhi}
    \Phi^{n+1}_i - \Phi^{n}_i = \frac{\Delta t}{2}( \Pi^{n}_i + \Pi^{n+1}_i) + J_0(\Pi)_i ~ \frac{\Delta t - 2 \Delta \tau}{2} + J_1(\Pi)_i ~ \frac{\Delta \tau}{2} ( \Delta \tau - \Delta t )
\end{equation}
where $\Delta \tau = \tau_i - t_n $ and $J_0(\Pi), J_1(\Pi)$ denote the jump discontinuities in $\Pi, \dot \Pi$ respectively. The subscript $i$ indicates that these terms are only applied at position $x_i$ if the particle worldline $\xi(t)$ crosses the $i$-th gridpoint at a time $\tau_i \in [t_n,t_{n+1}]$: $\xi(\tau_i)=x_i$.

When applied to Eq. \eqref{spacePi}, the discontinuous trapezium rule  \eqref{eq:discTrap} yields
\begin{multline}\label{discretePi}
    \Pi_i^{n+1} - \Pi_i^{n} = \frac{\Delta t}{2} \Bar{D}^{(2)}_{ij} (\Phi_j^{n} + \Phi_j^{n+1}) + \frac{\Delta t}{2} M_{ij} (\Pi_j^{n} + \Pi_j^{n+1}) + \frac{\Delta t}{2}( \Bar{s}_i(t_n) + \Bar{s}_i(t_{n+1}))\\
    + J_0(\partial_x^2\Phi)_i ~ \frac{\Delta t - 2 \Delta \tau}{2} + J_1(\partial_x^2\Phi)_i ~ \frac{\Delta \tau}{2} ( \Delta \tau - \Delta t ) + \llbracket \Pi \rrbracket_i
\end{multline}
where we have introduced the notation
\begin{equation}
    \llbracket \Pi \rrbracket_i = - \int_{t_n}^{t_{n+1}} \cos{t}~\delta(x_i - vt) ~ dt.
\end{equation}
Like the $J_k$ terms, this only  turns on if the particle worldline $\xi(t)=v t$ crosses $x_i$ at a time $\tau_i \in [t_n,t_{n+1}]$: $v \tau_i=x_i$.

This is now a linear system which can be solved for $\Phi^{n+1}$ and $\Pi^{n+1}$. When Eq. \eqref{discretePhi} is substituted into Eq. \eqref{discretePi}, we find
\begin{multline}
    \Pi_i^{n+1} - \Pi^{n}_i = \Delta t \Bar{D}^{(2)}_{ij} \Phi_j^{n} + \frac{\Delta t^2}{4} \Bar{D}^{(2)}_{ij}( \Pi^{n}_j + \Pi^{n+1}_j) + \frac{\Delta t}{2} M_{ij} ( \Pi^{n}_j + \Pi^{n+1}_j)\\
    + \frac{\Delta t}{2}( \Bar{s}_i(t_n) + \Bar{s}_i(t_{n+1})) + \frac{\Delta t}{2} \Bar{D}^{(2)}_{ij}\bigg( J_0(\Pi)_j ~ \frac{\Delta t - 2 \Delta \tau}{2} + J_1(\Pi)_j ~ \frac{\Delta \tau}{2} ( \Delta \tau - \Delta t ) \bigg)\\
    + J_0(\partial_x^2\Phi)_i ~ \frac{\Delta t - 2 \Delta \tau}{2} + J_1(\partial_x^2\Phi)_i ~ \frac{\Delta \tau}{2} ( \Delta \tau - \Delta t ) + \llbracket \Pi \rrbracket_i
\end{multline}
Some algebraic manipulations reveal that this scheme may be expressed as a matrix equation:
\begin{multline}
    \bigg( I_{ij} - \frac{\Delta t}{2} M_{ij} - \frac{\Delta t^2}{4} \Bar{D}_{ij}^{(2)} \bigg) \Pi_{j}^{n+1} = \bigg( I_{ij} + \frac{\Delta t}{2} M_{ij} + \frac{\Delta t^2}{4} \Bar{D}_{ij}^{(2)} \bigg) \Pi_{j}^{n} + \Delta t \Bar{D}^{(2)}_{ij} \Phi_j^{n}\\
    + \frac{\Delta t}{2}( \Bar{s}_i(t_n) + \Bar{s}_i(t_{n+1})) + \frac{\Delta t}{2} \Bar{D}^{(2)}_{ij}\bigg( J_0(\Pi)_j ~ \frac{\Delta t - 2 \Delta \tau}{2} + J_1(\Pi)_j ~ \frac{\Delta \tau}{2} ( \Delta \tau - \Delta t ) \bigg)\\
    + J_0(\partial_x^2\Phi)_i ~ \frac{\Delta t - 2 \Delta \tau}{2} + J_1(\partial_x^2\Phi)_i ~ \frac{\Delta \tau}{2} ( \Delta \tau - \Delta t ) + \llbracket \Pi \rrbracket_i
\end{multline}
All terms on the right side of this equation are quantities either known at $t_n$ or that can be computed in advance. Thus, solving for $\mathbf{\Pi}^{n+1}$ is tantamount to solving a linear system of the form $\mathbf{A}\cdot \mathbf{x} = \mathbf{b}$. The spatial discretization used determines the best way to proceed. If a finite difference or finite element scheme is used, then the matrix $\mathbf{A}=\mathbf{I}-\tfrac{\Delta t}{2} \mathbf{M}-\tfrac{\Delta t^2}{4}  \mathbf{\bar D}^2$ multiplying $\mathbf{\Pi}^{n+1}$  is sparse, in which case a direct row reduction method such as the Thomas algorithm scales linearly in $N$. If a pseudo-spectral method is used (as we do), then the matrix $\mathbf{A}$ multiplying $\mathbf{\Pi}^{n+1}$ is full and it is more efficient to invert it in advance and store it, since matrix-vector multiplication scales like $N^2$ while row-reducing a full matrix scales like $N^3$.
Once $\mathbf{\Pi}^{n+1}$ is known, it may be substituted into Eq. \eqref{discretePhi} to determine $\mathbf{\Phi}^{n+1}$. This is the time-stepping algorithm used to evolve the system.

As a final comment, we note that the condition number of $\mathbf{\bar D}^{(n)}$ (for the Chebyschev collocation method described in §2.2) is $\kappa(\mathbf{\bar D}^{(n)})\sim CN^{2n}$; for $N=100$ and $n=2$ this is approximately $2\times 10^6$. However, for the matrix $\mathbf{A}$ being inverted we find that, for $\Delta t =10^{-1}$,  $\kappa(\mathbf{A}) \sim 500$, while, for $\Delta t =10^{-5}$,  $\kappa(\mathbf{A}) \sim 1$. Thus, the fact that the  matrix $\mathbf{A}$ is well-conditioned as 
$\Delta t \rightarrow 0$ justifies the inversion and storage procedure outlined above.

\section*{References}

\bibliographystyle{siam}
\bibliography{library}

\begin{thebibliography}{10}

\bibitem{Adcock2010}
{\sc B.~Adcock}, {\em {Convergence acceleration of modified Fourier series in
  one or more dimensions}}, Math. Comput., 80 (2010), pp.~225--261.

\bibitem{Archibald2003}
{\sc R.~Archibald, K.~Chen, A.~Gelb, and R.~Renaut}, {\em {Improving tissue
  segmentation of human brain MRI through preprocessing by the Gegenbauer
  reconstruction method}}, Neuroimage, 20 (2003), pp.~489--502.

\bibitem{Archibald2002}
{\sc R.~Archibald and A.~Gelb}, {\em {A method to reduce the Gibbs ringing
  artifact in MRI scans while keeping tissue boundary integrity.}}, IEEE Trans.
  Med. Imaging, 21 (2002), pp.~305--19.

\bibitem{Archibald2004}
{\sc R.~Archibald, J.~Hu, A.~Gelb, and G.~Farin}, {\em {Improving the Accuracy
  of Volumetric Segmentation Using Pre-Processing Boundary Detection and Image
  Reconstruction}}, IEEE Trans. Image Process., 13 (2004), pp.~459--466.

\bibitem{Baltensperger2003}
{\sc R.~Baltensperger and M.~R. Trummer}, {\em {Spectral Differencing with a
  Twist}}, SIAM J. Sci. Comput., 24 (2003), pp.~1465--1487.

\bibitem{Barack2009}
{\sc L.~Barack}, {\em {Gravitational self-force in extreme mass-ratio
  inspirals}}, Class. Quantum Gravity, 26 (2009), p.~213001.

\bibitem{Barack:2000zq}
{\sc L.~Barack and L.~M. Burko}, {\em {Radiation reaction force on a particle
  plunging into a black hole}}, Phys. Rev. D, 62 (2000), p.~084040.

\bibitem{Giudice2015}
{\sc L.~Barack and P.~Giudice}, {\em Time-domain metric reconstruction for
  self-force applications}, Phys. Rev. D, 95 (2017), p.~104033.

\bibitem{Barack:2002ku}
{\sc L.~Barack and C.~O. Lousto}, {\em {Computing the gravitational selfforce
  on a compact object plunging into a Schwarzschild black hole}}, Phys. Rev. D,
  66 (2002), p.~061502.

\bibitem{Barack:2005nr}
\leavevmode\vrule height 2pt depth -1.6pt width 23pt, {\em {Perturbations of
  Schwarzschild black holes in the Lorenz gauge: Formulation and numerical
  implementation}}, Phys. Rev. D, 72 (2005), p.~104026.

\bibitem{Barkhudaryan2007}
{\sc A.~Barkhudaryan, R.~Barkhudaryan, and A.~Poghosyan}, {\em {Asymptotic
  behavior of Eckhoff’s method for Fourier series convergence acceleration}},
  Anal. Theory Appl., 23 (2007), pp.~228--242.

\bibitem{Batenkov2012}
{\sc D.~Batenkov and Y.~Yomdin}, {\em {Algebraic Fourier reconstruction of
  piecewise smooth functions}}, Math. Comput., 81 (2012), pp.~277--318.

\bibitem{Fornberg1998}
{\sc {Bengt Fornberg}}, {\em {A Practical Guide to Pseudospectral Methods}},
  Cambridge University Press, 1998.

\bibitem{Birkhoff1906}
{\sc G.~D. Birkhoff}, {\em {General mean value and remainder theorems with
  applications to mechanical differentiation and quadrature}}, Trans. Am. Math.
  Soc., 7 (1906), pp.~107--107.

\bibitem{Boyd2001}
{\sc J.~P. Boyd}, {\em {Chebyshev and Fourier Spectral Methods: Second Revised
  Edition (Dover Books on Mathematics)}}, Dover Publications, 2001.

\bibitem{Boyd2005}
\leavevmode\vrule height 2pt depth -1.6pt width 23pt, {\em {Trouble with
  Gegenbauer reconstruction for defeating Gibbs’ phenomenon: Runge phenomenon
  in the diagonal limit of Gegenbauer polynomial approximations}}, J. Comput.
  Phys., 204 (2005), pp.~253--264.

\bibitem{Burko2000}
{\sc L.~M. Burko}, {\em {Self-force on static charges in Schwarzschild
  spacetime}}, Class. Quantum Gravity, 17 (2000), pp.~227--250.

\bibitem{Canizares2009}
{\sc P.~Canizares and C.~Sopuerta}, {\em {Efficient pseudospectral method for
  the computation of the self-force on a charged particle: Circular geodesics
  around a Schwarzschild black hole}}, Phys. Rev. D, 79 (2009), p.~084020.

\bibitem{Canizares2011a}
{\sc P.~Canizares and C.~F. Sopuerta}, {\em {Tuning time-domain pseudospectral
  computations of the self-force on a charged scalar particle}}, Class. Quantum
  Gravity, 28 (2011), p.~134011.

\bibitem{Canizares2010c}
{\sc P.~Canizares, C.~F. Sopuerta, and J.~L. Jaramillo}, {\em {Pseudospectral
  collocation methods for the computation of the self-force on a charged
  particle: Generic orbits around a Schwarzschild black hole}}, Phys. Rev. D,
  82 (2010), p.~044023.

\bibitem{Costa2000}
{\sc B.~Costa and W.~S. Don}, {\em {On the computation of high order
  pseudospectral derivatives}}, Appl. Numer. Math., 33 (2000), pp.~151--159.

\bibitem{Csukas2019kcb}
{\sc K.~Csukas, I.~Racz, and G.~Z. Toth}, {\em {Numerical investigation of the
  dynamics of linear spin $s$ fields on a Kerr background: Late-time tails of
  spin $s = \pm 1, \pm 2$ fields}}, Phys. Rev. D, 100 (2019), p.~104025.

\bibitem{Dolan2013}
{\sc S.~R. Dolan}, {\em {Superradiant instabilities of rotating black holes in
  the time domain}}, Phys. Rev. D, 87 (2013), p.~124026.

\bibitem{Don1997}
{\sc W.~S. Don and A.~Solomonoff}, {\em {Accuracy Enhancement for Higher
  Derivatives using Chebyshev Collocation and a Mapping Technique}}, SIAM J.
  Sci. Comput., 18 (1997), pp.~1040--1055.

\bibitem{eckhoff1994}
{\sc K.~S. Eckhoff}, {\em {On discontinuous solutions of hyperbolic
  equations}}, Comput. Methods Appl. Mech. Engrg., 116 (1994), pp.~103--112.

\bibitem{eckhoff1995}
\leavevmode\vrule height 2pt depth -1.6pt width 23pt, {\em {Accurate
  reconstructions of functions of finite regularity from truncated Fourier
  series expansions}}, Math. Comp., 64 (1995), pp.~671--690.

\bibitem{eckhoff1997}
\leavevmode\vrule height 2pt depth -1.6pt width 23pt, {\em {On a high order
  numerical method for solving partial differential equations in complex
  geometries}}, J. Sci. Comput., 12 (1997), pp.~119--138.

\bibitem{eckhoff1998}
\leavevmode\vrule height 2pt depth -1.6pt width 23pt, {\em {On a high order
  numerical method for functions with singularities}}, Math. Comp., 67 (1998),
  pp.~1063--1087.

\bibitem{eckhoff1996}
{\sc K.~S. Eckhoff and J.~H. Rolfsnes}, {\em {On nonsmooth solutions of linear
  hyperbolic systems}}, J. Comput. Phys., 125 (1996), pp.~1--15.

\bibitem{Field2009}
{\sc S.~E. Field, J.~S. Hesthaven, and S.~R. Lau}, {\em {Discontinuous Galerkin
  method for computing gravitational waveforms from extreme mass ratio
  binaries}}, Class. Quantum Gravity, 26 (2009), p.~165010.

\bibitem{Field:2010xn}
{\sc S.~E. Field, J.~S. Hesthaven, and S.~R. Lau}, {\em {Persistent junk
  solutions in time-domain modeling of extreme mass ratio binaries}}, Phys.
  Rev. D, 81 (2010), p.~124030.

\bibitem{Fornberg1988}
{\sc B.~Fornberg}, {\em {Generation of finite difference formulas on
  arbitrarily spaced grids}}, Math. Comput., 51 (1988), pp.~699--699.

\bibitem{Fornberg2006}
{\sc B.~{}Fornberg}, {\em {A Pseudospectral Fictitious Point Method for High
  Order Initial Boundary Value Problems}}, SIAM J. Sci. Comput., 28 (2006),
  pp.~1716--1729.

\bibitem{GottliebDavidSigal2003}
{\sc D.~Gottlieb and S.~Gottlieb}, {\em {Spectral Methods for Discontinuous
  Problems}}, in NA03 Dundee, 2003.

\bibitem{Gottlieb1994}
{\sc D.~Gottlieb and C.-W. Shu}, {\em {Resolution properties of the Fourier
  method for discontinuous waves}}, Comput. Methods Appl. Mech. Eng., 116
  (1994), pp.~27--37.

\bibitem{Gottlieb1995}
\leavevmode\vrule height 2pt depth -1.6pt width 23pt, {\em {On the Gibbs
  phenomenon IV. Recovering exponential accuracy in a subinterval from a
  Gegenbauer partial sum of a piecewise analytic function}}, Math. Comput., 64
  (1995), pp.~1081--1081.

\bibitem{Gottlieb1995a}
\leavevmode\vrule height 2pt depth -1.6pt width 23pt, {\em {On the Gibbs
  phenomenon V: recovering exponential accuracy from collocation point values
  of a piecewise analytic function}}, Numer. Math., 71 (1995), pp.~511--526.

\bibitem{Gottlieb1996}
\leavevmode\vrule height 2pt depth -1.6pt width 23pt, {\em {On the Gibbs
  Phenomenon III: Recovering Exponential Accuracy in a Sub-Interval From a
  Spectral Partial Sum of a Pecewise Analytic Function}}, SIAM J. Numer. Anal.,
  33 (1996), pp.~280--290.

\bibitem{Gottlieb1997}
\leavevmode\vrule height 2pt depth -1.6pt width 23pt, {\em {On the Gibbs
  Phenomenon and Its Resolution}}, SIAM Rev., 39 (1997), pp.~644--668.

\bibitem{Gottlieb1992}
{\sc D.~Gottlieb, C.-W. Shu, A.~Solomonoff, and H.~Vandeven}, {\em {On the
  Gibbs phenomenon I: recovering exponential accuracy from the Fourier partial
  sum of a nonperiodic analytic function}}, J. Comput. Appl. Math., 43 (1992),
  pp.~81--98.

\bibitem{Haas:2007kz}
{\sc R.~Haas}, {\em {Scalar self-force on eccentric geodesics in Schwarzschild
  spacetime: A Time-domain computation}}, Phys. Rev. D, 75 (2007), p.~124011.

\bibitem{Haas:2011np}
\leavevmode\vrule height 2pt depth -1.6pt width 23pt, {\em {Time domain
  calculation of the electromagnetic self-force on eccentric geodesics in
  Schwarzschild spacetime}},  (2011).

\bibitem{2010Lubina}
{\sc E.~{Hairer}, C.~{Lubich}, and G.~{Wanner}}, {\em {Geometric Numerical
  Integration: Structure-Preserving Algorithms for Ordinary Differential
  Equations}}, Springer,  (2010).

\bibitem{Harms2013}
{\sc E.~Harms, S.~Bernuzzi, and B.~Bruegmann}, {\em {Numerical solution of the
  2+1 Teukolsky equation on a hyperboloidal and horizon penetrating foliation
  of Kerr and application to late-time decays}},  (2013).

\bibitem{Harten1989}
{\sc A.~Harten}, {\em {ENO schemes with subcell resolution}}, J. Comput. Phys.,
  83 (1989), pp.~148--184.

\bibitem{Heffernan:2017cad}
{\sc A.~Heffernan, A.~C. Ottewill, N.~Warburton, B.~Wardell, and P.~Diener},
  {\em {Accelerated motion and the self-force in Schwarzschild spacetime}},
  Class. Quant. Grav., 35 (2018), p.~194001.

\bibitem{hermite1878}
{\sc C.~Hermite}, {\em {Formule d'interpolation de Lagrange}}, J. reine angew.
  Math, 84 (1878), pp.~70--79.

\bibitem{Jaramillo2011}
{\sc J.~L. Jaramillo, C.~F. Sopuerta, and P.~Canizares}, {\em {Are time-domain
  self-force calculations contaminated by Jost solutions?}}, Phys. Rev. D, 83
  (2011), p.~061503.

\bibitem{Jung2004}
{\sc J.-H. Jung and B.~D. Shizgal}, {\em {Generalization of the inverse
  polynomial reconstruction method in the resolution of the Gibbs phenomenon}},
  J. Comput. Appl. Math., 172 (2004), pp.~131--151.

\bibitem{kantorovichkrylov1958}
{\sc L.~V. Kantorovich and V.~I. Krylov}, {\em {Approximate methods of higher
  analysis}},  (1958).

\bibitem{2020Wolfram}
{\sc R.~Knapp}, {\em {Wolfram Language Documentation - The Numerical Method of
  Lines}}.
\newblock
  \url{https://reference.wolfram.com/language/tutorial/NDSolveMethodOfLines.ht%
ml}.
\newblock [Accessed December 25, 2020].

\bibitem{knapp2008method}
{\sc R.~Knapp}, {\em {A Method of Lines Framework in Mathematica}}, JNAIAM, 3
  (2008), pp.~43--59.

\bibitem{Kosloff1993}
{\sc D.~Kosloff and H.~Tal-Ezer}, {\em {A Modified Chebyshev Pseudospectral
  Method with an O(N-1) Time Step Restriction}}, J. Comput. Phys., 104 (1993),
  pp.~457--469.

\bibitem{Krylov1907}
{\sc A.~N. Krylov}, {\em {On approximate calculations, Lectures delivered in
  1906}}, Tipolitography of Birkenfeld, St. Petersburg, 1907.

\bibitem{Lanczos1966}
{\sc C.~Lanczos}, {\em {Discourse on Fourier series}}, Oliver \& Boyd, 1966.

\bibitem{Lipman2009}
{\sc Y.~Lipman and D.~Levin}, {\em {Approximating piecewise-smooth functions}},
  IMA J. Numer. Anal., 30 (2009), pp.~1159--1183.

\bibitem{Lousto:1999za}
{\sc C.~O. Lousto}, {\em {Pragmatic approach to gravitational radiation
  reaction in binary black holes}}, Phys. Rev. Lett., 84 (2000),
  pp.~5251--5254.

\bibitem{Lousto:1997wf}
{\sc C.~O. Lousto and R.~H. Price}, {\em {Understanding initial data for black
  hole collisions}}, Phys. Rev. D, 56 (1997), pp.~6439--6457.

\bibitem{2019arXiv190109967M}
{\sc C.~M. {Markakis}, M.~F. {O'Boyle}, D.~{Glennon}, K.~{Tran}, P.~{Brubeck},
  R.~{Haas}, H.-Y. {Schive}, and K.~{Ury{\={u}}}}, {\em {Time-symmetry,
  symplecticity and stability of Euler-Maclaurin and Lanczos-Dyche
  integration}}, arXiv e-prints,  (2019), p.~arXiv:1901.09967.

\bibitem{Martel:2001yf}
{\sc K.~Martel and E.~Poisson}, {\em {A One parameter family of time symmetric
  initial data for the radial infall of a particle into a Schwarzschild black
  hole}}, Phys. Rev. D, 66 (2002), p.~084001.

\bibitem{eckhoff2002}
{\sc O.~F. N\ae~ss and K.~S. Eckhoff}, {\em {A modified Fourier-Galerkin method
  for the Poisson and Helmholtz equations}}, in J. Sci. Comput., vol.~17, 2002,
  pp.~529--539.

\bibitem{2019JCoPh.390..527P}
{\sc J.~{Piotrowska}, J.~M. {Miller}, and E.~{Schnetter}}, {\em {Spectral
  methods in the presence of discontinuities}}, Journal of Computational
  Physics, 390 (2019), pp.~527--547.

\bibitem{Poghosyan2010}
{\sc A.~Poghosyan}, {\em {Asymptotic behavior of the Eckhoff method for
  convergence acceleration of trigonometric interpolation}}, Anal. Theory
  Appl., 26 (2010), pp.~236--260.

\bibitem{Poghosyan2012}
\leavevmode\vrule height 2pt depth -1.6pt width 23pt, {\em {On an
  autocorrection phenomenon of the Eckhoff interpolation}}, Aust. J. Math.
  Anal. Appl., 9 (2012), pp.~1 -- 31.

\bibitem{Poisson2011}
{\sc E.~Poisson, A.~Pound, and I.~Vega}, {\em {The Motion of Point Particles in
  Curved Spacetime}}, Living Rev. Relativ., 14 (2011).

\bibitem{Racz2011}
{\sc I.~R\'{a}cz and G.~Z. T\'{o}th}, {\em {Numerical investigation of the
  late-time Kerr tails}}, Class. Quantum Gravity, 28 (2011), p.~195003.

\bibitem{Regge1957}
{\sc T.~Regge and J.~A. Wheeler}, {\em {Stability of a Schwarzschild
  Singularity}}, Phys. Rev., 108 (1957), pp.~1063--1069.

\bibitem{Sadiq2011}
{\sc B.~Sadiq and D.~Viswanath}, {\em {Finite Difference Weights, Spectral
  Differentiation, and Superconvergence}}, Math. Comput.,  (2011).

\bibitem{Schiesser1991}
{\sc W.~E. Schiesser}, {\em {The Numerical Method of Lines: Integration of
  Partial Differential Equations}}, Academic Press, 1991.

\bibitem{Shu1995}
{\sc C.-W. Shu and P.~S. Wong}, {\em {A note on the accuracy of spectral method
  applied to nonlinear conservation laws}}, J. Sci. Comput., 10 (1995),
  pp.~357--369.

\bibitem{Sundararajan2007}
{\sc P.~Sundararajan, G.~Khanna, and S.~Hughes}, {\em {Towards adiabatic
  waveforms for inspiral into Kerr black holes: A new model of the source for
  the time domain perturbation equation}}, Phys. Rev. D, 76 (2007), p.~104005.

\bibitem{Sundararajan2008}
{\sc P.~Sundararajan, G.~Khanna, S.~Hughes, and S.~Drasco}, {\em {Towards
  adiabatic waveforms for inspiral into Kerr black holes. II. Dynamical sources
  and generic orbits}}, Phys. Rev. D, 78 (2008), p.~024022.

\bibitem{Teukolsky1973}
{\sc S.~A. Teukolsky}, {\em {Perturbations of a Rotating Black Hole. I.
  Fundamental Equations for Gravitational, Electromagnetic, and Neutrino-Field
  Perturbations}}, Astrophys. J., 185 (1973), p.~635.

\bibitem{Welfert1997}
{\sc B.~D. Welfert}, {\em {Generation of Pseudospectral Differentiation
  Matrices I}}, SIAM J. Numer. Anal., 34 (1997), pp.~1640--1657.

\bibitem{2011JCoPh.230.2286Z}
{\sc A.~{Zengino{\u{g}}lu}}, {\em {Hyperboloidal layers for hyperbolic
  equations on unbounded domains}}, Journal of Computational Physics, 230
  (2011), pp.~2286--2302.

\bibitem{Zerilli1970a}
{\sc F.~J. Zerilli}, {\em {Effective Potential for Even-Parity Regge-Wheeler
  Gravitational Perturbation Equations}}, Phys. Rev. Lett., 24 (1970),
  pp.~737--738.

\bibitem{Zerilli1970}
\leavevmode\vrule height 2pt depth -1.6pt width 23pt, {\em {Gravitational Field
  of a Particle Falling in a Schwarzschild Geometry Analyzed in Tensor
  Harmonics}}, Phys. Rev. D, 2 (1970), pp.~2141--2160.

\end{thebibliography}

\end{document}